\documentclass[11pt,a4paper]{article}
\usepackage{amsmath}
\usepackage{amsmath,amscd}
\usepackage{amssymb,amsfonts,amsmath,amscd}
\usepackage[margin=0.5in,footskip=0.25in]{geometry}
\usepackage{array}
\usepackage{diagbox}
\usepackage{makecell}
\usepackage{hhline}
\usepackage{tikz}
\usetikzlibrary{arrows}
 \usetikzlibrary{calc}
 \usepackage{graphicx}

\tikzset{me/.style={to path={
\pgfextra{%
 \pgfmathsetmacro{\startf}{-(#1-1)/2}
 \pgfmathsetmacro{\endf}{-\startf}
 \pgfmathsetmacro{\stepf}{\startf+1}}
 \ifnum 1=#1 -- (\tikztotarget)  \else
     let \p{mid}=($(\tikztostart)!0.5!(\tikztotarget)$)
         in
\foreach \i in {\startf,\stepf,...,\endf}
    {%
     (\tikztostart) .. controls ($ (\p{mid})!\i*10pt!90:(\tikztotarget) $) .. (\tikztotarget)
      }
      \fi
     \tikztonodes
}}}

\usepackage{fancybox}
\usepackage{pgf,pgfarrows,pgfnodes,pgfautomata,pgfheaps,pgfshade}
\usepackage[all,graph,frame]{xy}
\usepackage{appendix}
\usepackage{enumerate}

\newcommand{\raf}[1]{(\ref{#1})}

\newcommand{\cC}{{\cal C}}
\newcommand{\cD}{{\cal D}}
\newcommand{\cE}{{\cal E}}

\newcommand{\cH}{{\cal H}}

\def\R{\mathbb{R}}

\def\proof{{\em Proof}}
\def\qed{{$\Box$}}

\newtheorem{theorem}{Theorem}

\newtheorem{proposition}{Proposition}
\newtheorem{remark}{Remark}

\begin{document}


\title{ Monotone bargaining is Nash-solvable}
\author{Vladimir Gurvich\footnote{ National Research University, Higher School of Economics (HSE), Moscow Russia
\textit{vgurvich@hse.ru; vladimir.gurvich@gmal.com}} and
Gleb Koshevoy\footnote{CEMI, Russian Academy of Science of Moscow
\textit{e-mail: koshevoyga@gmail.com}}}
\maketitle



\begin{abstract}
Given two finite ordered sets
$A = \{a_1, \ldots, a_m\}$  and  $B = \{b_1, \ldots, b_n\}$,
introduce the set of  $m n$  outcomes of the game
$O = \{(a, b) \mid a \in A, b \in B\} =
\{(a_i, b_j) \mid  i \in I = \{1, \ldots, m\},  j \in J = \{1, \ldots, n\}$.
Two players, Alice and Bob, have the sets of strategies $X$  and  $Y$
that consist of all monotone non-decreasing mappings
$x: A \rightarrow B$  and  $y: B \rightarrow A$, respectively.
It is easily seen that each pair $(x,y) \in X \times Y$
produces at least one {\em deal}, that is,
an outcome  $(a,b) \in O$  such that $x(a) = b$  and  $y(b) = a$.
Denote by  $G(x,y) \subseteq O$ the set of all such deals related to $(x,y)$.
The obtained mapping  $G = G_{m,n}: X \times Y \rightarrow 2^O$
is a {\em game correspondence}.
Choose an arbitrary deal $g(x,y) \in G(x,y)$
to obtain a mapping  $g : X \times Y \rightarrow O$,
which is a {\em game form}. We use notation  $g \in G$  and show that 
each such game form is tight and, hence, Nash-solvable, that is,
for any pair  $u = (u_A, u_B)$  of utility functions
$u_A : O \rightarrow \R$  of Alice and  $u_B: O \rightarrow \R$ of Bob,
the obtained monotone bargaining game  $(g, u)$  has at least one
Nash equilibrium in pure strategies.
Moreover, the same equilibrium can be selected for all $g \in G$.
We also obtain an efficient algorithm that determines such an equilibrium
in time linear in $m n$, although the numbers of strategies
$|X| = \binom{m+n-1}{m}$ and $|Y| = \binom{m+n-1}{n}$  are exponential in $m n$.
Our results show that,  somewhat surprising, the players have no need to hide
or randomize their bargaining strategies, even in the zero-sum case.
\end{abstract}


\section{Introduction}
\label{s0}
\subsection{Main results}
\label{ss00}
Two players, Alice and Bob, possess items
$A = \{a_1, \ldots, a_m\}$  and  $B = \{b_1, \ldots, b_n\}$, respectively.
Both sets are ordered: the smaller is the index,
the less valuable is the corresponding item.
Both players know both orders.
A pair $(a, b)$  with $a \in A$ and $b \in B$  is called an outcome.

Alice chooses her strategy, which is a mappings  $x: A \rightarrow B$,
showing that she is ready to exchange  $a$  for  $x(a)$   for any  $a \in A$.
Naturally, such a mapping  $x$  is assumed to be {\em monotone non-decreasing}:
$x(a) \geq  x(a')$  whenever  $a > a'$.
Similarly, Bob chooses his strategy, which is
a monotone non-decreasing mapping  $y: B \rightarrow A$.
It is not difficult to compute the cardinalities
of the sets of strategies and outcomes:

\begin{equation}
\label{eq1}
|X| = \binom {m+n-1}{m}, \; |Y| = \binom {m+n-1}{n}; \;\; |O| = mn.
\end{equation}

Let us fix a pair of strategies  $(x,y)$.
An outcome  $(a,b) \in O$ is called  a {\em deal} if  $x(a) = b$  and  $y(b) = a$.
Denote by  $G(x,y) \subseteq O$  the set of all deals.
We will show that set  $G(x,y)$  is not empty;
yet, it may contain several deals.
The obtained mapping  $G : X \times Y \rightarrow 2^O$
is called a (monotone bargaining) {\em game correspondence}.
Note that  $G = G_{m,n}$  is uniquely defined by  $m$  and  $n$.
For each $(x,y)$, let us fix an arbitrary deal  $g(x,y)$  from  $G(x,y)$.
The obtained a mapping  $g : X \times Y \rightarrow O$,
is called a (monotone bargaining) {\em game form}.

An arbitrary (not necessarily monotone bargaining)
two-person game form  $g$  is called {\em Nash-solvable} (NS)
if for any pair  $u = (u_A, u_B)$  of utility functions
$u_A : O \rightarrow \R$  of Alice and  $u_B: O \rightarrow \R$ of Bob
the obtained game  $(g, u)$  has at least one
Nash equilibrium (NE) in pure strategies.
A game correspondence  $G$  is called NS  if every
its  game form  $g \in G$  is  NS.

We will apply some general criteria of Nash-solvability
(based on the concept of Boolean duality and
developed for the two-person game forms and correspondences
in \cite{Gur75}, \cite{Gur88}, \cite{LRMDS17})
to the monotone bargaining case considered in this paper.
These criteria imply that $G = G_{m,n}$  is NS  for all  $m$ and $n$.
Moreover, given $m$ and $n$, the same NE
$(x,y)$  exists for all $g \in G$,
and this equilibrium is simple, that is, $G(x,y)$  contains a unique deal.
We also construct an efficient algorithm that determines such an equilibrium
in time polynomial in $m n$, although,
by \raf{eq1}, the sizes of  $X$ and $Y$  are exponential in $m n$.

\medskip

In our model, Alice and  Bob are restricted
to their monotone non-decreasing strategies.
Let us replace them by the monotone non-increasing ones.
Obviously, this case is equivalent with the original one:
it is enough just to inverse the enumeration
of either  $A$, or  $B$, but not both.
However, if we restrict Alice
by her monotone non-decreasing and Bob by his
monotone non-increasing) strategies, or vice versa,
then  $G(x,y)$  may become empty for some  $x$  and  $y$.

\medskip

In particular, our results imply that
the players have no need to keep
their bargaining strategies secret and choose them randomly,
even in the zero-sum case.
This is unusual, since in most bargaining models players do
have to hide and randomize their strategies
\cite{DS91}, \cite{Nas50}, \cite{Nas53}, \cite{Rai53}, \cite{Rap70}, \cite{FS82}.

\smallskip

The considered monotone bargaining games
can be viewed as positional games on complete bipartite digraphs,
in which, however, both players are restricted to their monotone non-decreasing
(or equivalently, non-increasing) strategies.
Restricting the set of strategies can be viewed as
a generalization of the concept of imperfect information \cite[Section 6]{Gur78}.
Nevertheless, we show that
the monotone bargaining games can be solved in pure strategies.
It would be interesting to find other classes
of positional games with imperfect information
NS in pure strategies.


\medskip

The paper is organized as follows.
In Section \ref{s1} we derive some combinatorial properties
of a monotone bargaining game correspondence $G = G_{m,n}$  and then
explain their role and meaning, after recalling some
general definitions and results related to
the two-person normal form games in Section \ref{s2}.
In particular, we show that  $G$,
as well as all game forms   $g \in G$,  are tight.
Moreover, they have the same pair of dual hypergraphs
(or the same self-dual effectivity functions, in terminology of
\cite{Gur97} and \cite{Gur97a}).
It is known that tightness and Nash-solvability are
equivalent properties of the two-person game forms;
see  \cite{Gur75}, \cite{Gur88}, and also \cite{LRMDS17}.
In Section \ref{ss23} we recall
a constructive proof of this fact that is based on an algorithm
finding a NE in any game  $(g,u)$  in time linear
in the size of the set $O$ of its outcomes, whenever $g$ is tight.
Recall that  $|O| = mn$  for the monotone bargaining game forms.
In section \ref{s3} we recall other classes
of two-person tight game forms that appear from: 
positional structures, Section \ref{ss30}; 
self-dual hypergraphs, Section \ref{ss31});   
symmetric dual hypergraphs and  veto voting, Section \ref{ss32};  
topological structure, where tightness results
from Jordan's curve theorem, Section \ref{ss33}.
Finally, in Section \ref{s4} we show
that a natural extension of our results to the case of three players fails.

\subsection{Examples}
Three game correspondences $G = G_{m,n}$  are shown
in Figures 1,\, 2, and  3
for ($m=n=2$), ($m=3, n=2$), and ($m=n=3$) respectively.
The corresponding tables are of sizes  $3 \times 3$, $4 \times 6$, and $10 \times 10$,
respectively, and contain $4, 6,$ and $9$  outcomes, in accordance with \raf{eq1}.

Each of these three figures consists of two tables.
In the first one, items  $a \in A$  and  $b \in B$  are
represented as the vertices  of a bipartite digraph,
denoted by white and black discs that are always located
on the left and on the right, respectively.
The values of items increase with their indices  $i \in I$  and  $j \in J$,
the higher is a disc, the smaller is its index and
the cheaper is the corresponding item.
The strategies  (monotone non-decreasing mappings) $x$  and  $y$
are denoted by directed arcs, which
may have common ends, but no other intersections.
Indeed, they cannot have a common beginning, because
$x : A \rightarrow B$  and   $y : B \rightarrow A$  are mappings, and
since both mappings are monotone non-decreasing, the arcs cannot cross.
The strategies and deals are shown
on the corresponding bipartite graphs;
To save space, for every  $i \in I$  and  $j \in J$
we denote the deal $o_{ij} = (a_i, b_j)$
by a non-directed arc between $a_i$ and $b_j$.
In the second table the strategies and deals
are specified simply by
the corresponding pairs of indices  $(i,j)$.

Rows and columns of both tables are labeled
by the strategies  $x \in X$  and $y \in Y$,  respectively.
The set of deals  $G(x,y$  is shown for each entry  $(x,y)$.
These sets are always non-empty
(see Proposition \ref{p1} below) but may contain several deals.
For example, in Figure 3, three deals $\{a_i, b_i\} \mid i = 1,2,3\}$ are assigned
to the pair  $(x,y)$  with  $x(a_i) = b_i$  and  $y(b_i) = a_i$  for  $i = 1,2,3$;
and one or two deals are assigned to any other pair  $(x,y)$.

\newpage



\vskip1cm

Figure 1. Monotone bargaining game correspondence $G(2,2).$

\newpage

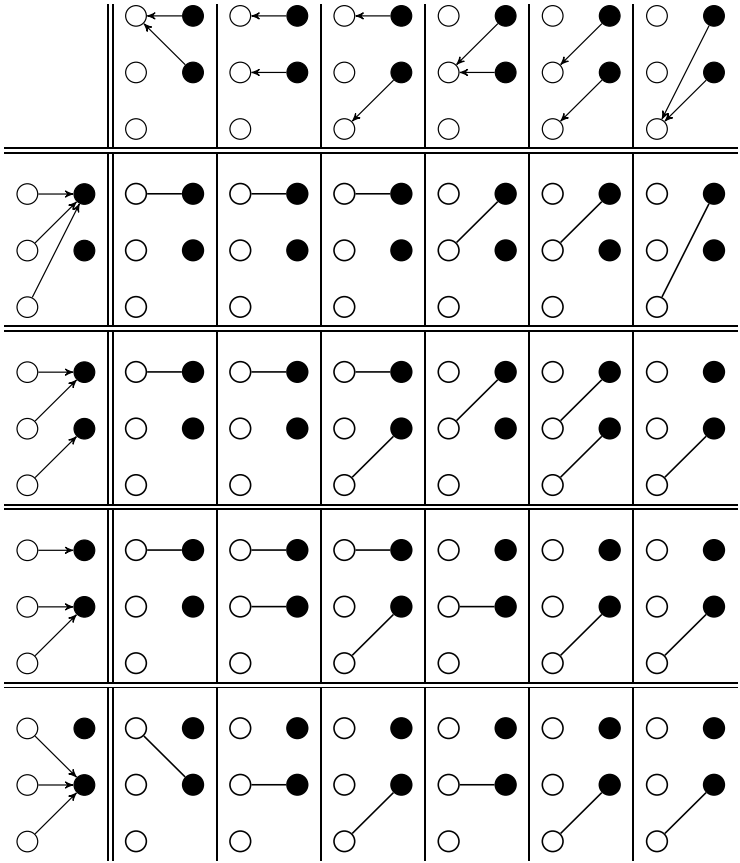
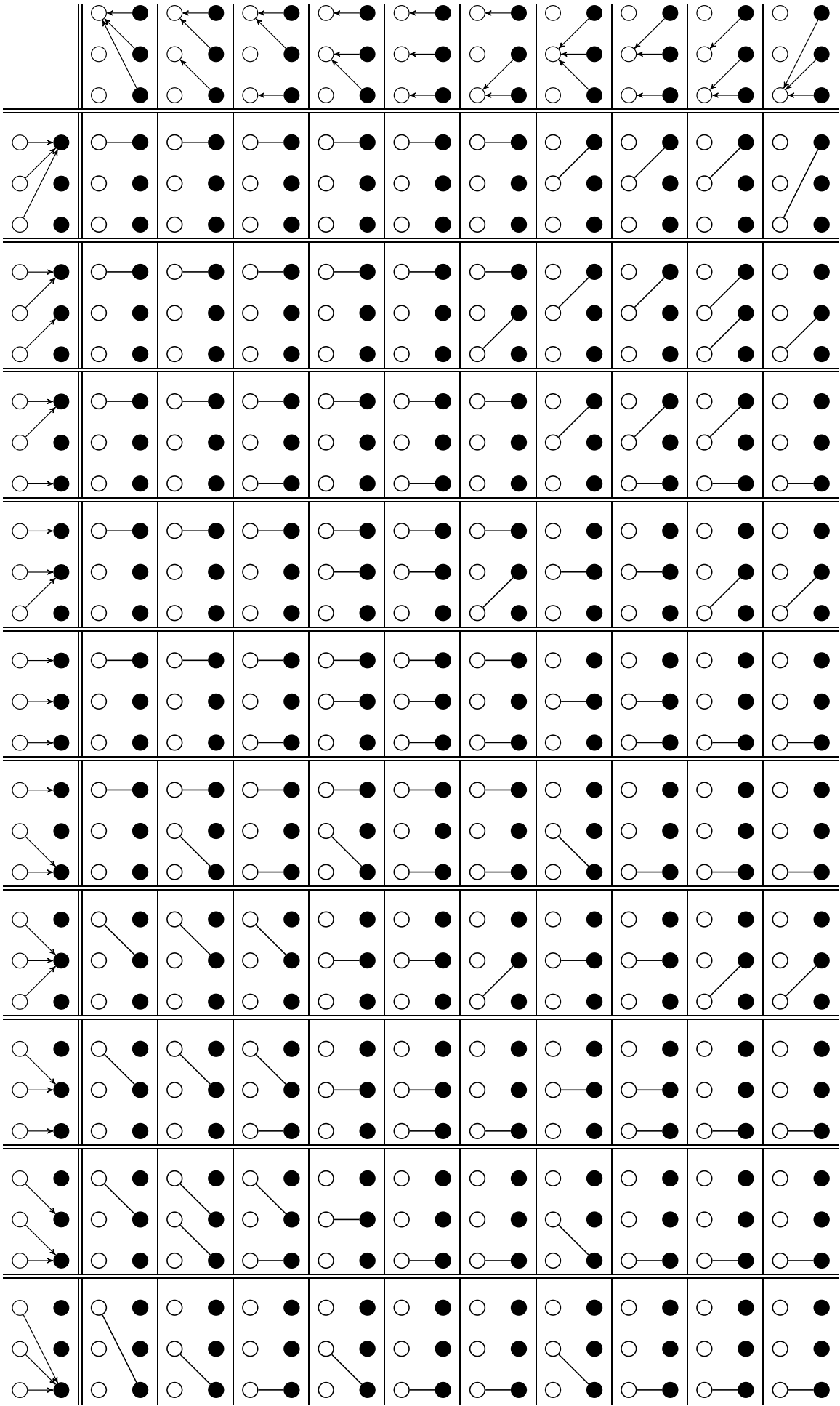
\scalebox{0.75}{
}

\vskip 2cm

{
\newcommand{\cl}[1]{\begin{tabular}{c c}#1 \end{tabular}}
\begin{tabular}{c||c|c|c|c|c|c}
& \cl{1 & 1 \\ 1 & 2} & \cl{1&1\\2&2} & \cl{1&1\\3&2 }& \cl{2&1\\2&2}& \cl{2&1\\3&2}& \cl{3&1\\3&2}
\\ \hhline{=======} & & & & & & \\[0.1em]
 \cl{ 1&1\\2&1\\3&2}& \cl{ 1&1}& \cl{ 1&1}& \cl{ 1&1}& \cl{ 2&1}& \cl{ 2&1}& \cl{ 3&1}
\\ \hhline{=======} & & & & & & \\[0.1em]
 \cl{ 1&1\\2&2\\3&2}& \cl{1&1 }& \cl{ 1&1}& \cl{ 1&1\\3&2}& \cl{ 2&1}& \cl{ 2&1\\3&2}& \cl{ 3&2}
\\ \hhline{=======} & & & & & & \\[0.1em]
 \cl{1&1\\2&2\\3&2 }& \cl{ 1&1}& \cl{ 1&1\\2&2}& \cl{ 1&1\\3&2}& \cl{2&2 }& \cl{ 3&2}& \cl{ 2&2}
 \\ \hhline{=======} & & & & & & \\[0.1em]
  \cl{ 1&2\\2&2\\3&2}& \cl{ 1&2}& \cl{ 2&2}& \cl{ 3&2}& \cl{ 2&2}& \cl{ 3&2}& \cl{ 3&2}
\end{tabular}
}

\bigskip

Figure 2. Monotone bargaining game correspondence $G(2, 3)$. 

\newpage

\scalebox{0.75}{
}

\newpage

{
\newcommand{\cl}[1]{\begin{tabular}{c c}#1 \end{tabular}}
\begin{tabular}{c||c|c|c|c|c|c|c|c|c|c}
& \cl{1&1\\1&2\\1&3} &\cl{1&1\\1&2\\2&3}&\cl{1&1\\1&2\\3&3}&\cl{1&1\\2&2\\2&3}&\cl{1&1\\2&2\\3&3}&\cl{1&1\\3&2\\3&3}&\cl{2&1\\2&2\\2&3}&\cl{2&1\\2&2\\3&3}&\cl{2&1\\3&2\\3&3}&\cl{3&1\\3&2\\3&3}
\\ \hhline{===========} & & & & & & & & & &\\[0.1em]
\cl{1&1\\2&1\\3&1}&\cl{1&1}&\cl{1&1}&\cl{1&1}&\cl{1&1}&\cl{1&1}&\cl{1&1}&\cl{2&1}&\cl{2&1}&\cl{2&1}&\cl{3&1}
\\ \hhline{===========} & & & & & & & & & &\\[0.1em]
\cl{1&1\\2&1\\3&2}&\cl{1&1}&\cl{1&1}&\cl{1&1}&\cl{1&1}&\cl{1&1}&\cl{1&1\\3&2}&\cl{2&1}&\cl{2&1}&\cl{2&1\\3&2}&\cl{3&3}
\\ \hhline{===========} & & & & & & & & & &\\[0.1em]
\cl{1&1\\2&1\\3&3}&\cl{1&1}&\cl{1&1}&\cl{1&1\\3&3}&\cl{1&1}&\cl{1&1\\3&3}&\cl{1&1}&\cl{2&1}&\cl{2&1\\3&3}&\cl{2&1\\3&3}&\cl{3&3}
\\ \hhline{===========} & & & & & & & & & &\\
\end{tabular}
}


\medskip

Figure 3. Monotone bargaining game correspondence $G(3, 3)$.

\bigskip 

For each of the three examples let us verify also the
following combinatorial properties, which will play an important role too:

\begin{itemize}
\item{(r)} For any deal  $(a,b) \in G(x, y)$
there exist  $x' \in X$  and  $y' \in Y$  such that
sets  $G(x, y')$  and  $G(x', y)$ contain only this deal  $(a,b)$.
\item{}In each row  $x \in X$ (respectively, column $y  \in Y$)
choose an arbitrary deal  $(a,b) \in G(x, y)$  for some  $y \in Y$
(respectively, for some $x  \in X$). Then,
\item{(t)} The obtained two selections always have a deal in common.
\item{(t$'$ $t''$)}
There exists a column $y_0$  (respectively, a row $x_0$)
such that the set of deals in it is a subset of
the above row-selection (respectively, column-selection).
\end{itemize}

Let us consider, for example, the $3 \times 3$ table in Figure 1.
The central entry $G(x_2,y_2)$ contains two deals
$(a_1,b_1)$  and  $(a_2,b_2)$, but the second row and the second column
each contains $(a_1,b_1)$  and  $(a_2,b_2)$, separately,
in  $G(x_1, y_2), G(x_2, y_1)$  and in $G(x_2, y_3), G(x_3, y_2)$.
Thus,  $(r)$  holds for  this table.
In the same table, the side diagonal is
a row selection that provides the set of deals
$(a_2, b_1), (a_1, b_1), (a_1, b_2)\}$, but the first
column provides $(a_1, b_1), (a_1, b_2)\}$, which is its subset.
Similarly, in the $10 \times 10$ table in Figure 3
the entry $G(x_5,y_5)$ contains three deals
$(a_1,b_1), (a_2,b_2)$,  and  $(a_3,b_3)$, but
each of these three deals appears in the fifth row and column
of this table separately, in agreement with  $(r)$.

In Sections \ref{s2} and \ref{s3}  we will recall
that properties $(t'), (t'')$ and $(t)$ are equivalent and imply $(r)$.
Due to $(r)$,  all  game forms  $g \in G$  have the same pair of hypergraphs
(or in other words, the same effectivity functions).
Furthermore, $(t'), (t'')$  and  $(t)$ mean that these hypergraphs are dual
(respectively, the effectivity function is self-dual)
implying that all game forms  $g \in G$  are tight and, hence,
Nash-solvable, in accordance with the criteria of \cite{Gur75} and \cite{Gur88}.

\subsection{Interpretation}
Let Alice and Bob be dealers possessing the sets of objects
$A$  and  $B$, respectively, and
a deal  $(a,b) \in A \times B$  means that
they exchange  $a$  and  $b$.
They may be art-dealers, car dealers; or one
of them may be just a buyer with a discrete budget.

It would be natural to order sets  $A$ and $B$
in accordance with the values of the objects.
This fully justifies the requirement of monotonicity
for the strategies  $x \in X$  and  $y \in Y$.
Yet, then, it is not natural to assume that
the utility functions $u_A$   and  $u_B$  are arbitrary.
Indeed, for a pair of deals  $(a,b)$ and  $(a',b')$
such that $a \leq a', b' \leq b$  and at least one
of these inequalities is strict, it would be natural to assume that
Alice strictly prefers  $(a,b)$  to $(a',b')$, while Bob vice versa.
Under such assumption, two constant strategies
$x : a \rightarrow b_n$   and   $y : b \rightarrow a_m$
form a trivial NE with the outcome  $(a_m,  b_n)$.
Although, we prove that an equilibrium exists for any  $u_A$  and $u_B$
(which is non-trivial) but why to consider any?

However, it becomes natural if the values of the deals
can be corrected by side-payments.
Then, let us also require that the total value
of each collection, $A$  and  $B$, should not be reduced.
These total values are estimated subjectively
by Alice and Bob (or by some experts) and
the side-payments are not taken into account.
Such a requirement would suffice to eliminate all strategies
except for the monotone non-decreasing ones.

Alternatively,  $A$ and $B$  are ordered in accordance with
some other attribute of the objects, distinct from their commercial value.
Such an attribute must be important enough
to justify the requirement of monotonicity of the strategies.
For example, it could be the shelf life of commodities,
or the age of pieces of art.
(The age is also a reasonable attribute
in case of the matching problems, when
$A$  and  $B$   are brides and grooms, and a deal is a marriage.)
We leave to the reader suggesting more practical cases
when  $x$  and  $y$  must be  monotone non-decreasing, while
$u$  can be arbitrary.

\section{Main properties of monotone bargaining}
\label{s1}

To any pair of strategies
$x : A \rightarrow B$  and  $y : B \rightarrow A$
(not necessarily  monotone non-decreasing)
let us assign a digraph  $\Gamma = \Gamma(x,y)$
on the vertex-set  $A \cup B$  as follows:  $(a,b)$
(respectively, $(b,a)$  is an arc of $\Gamma(x,y)$  whenever  $x(a) = b$
(respectively, $y(b) = a$).

Some visualization helps.
Let us consider an embedding of $\Gamma(x,y)$ into the plane,
as in Figures 1,\, 2, and 3, where many
($3 \times 3 + 4 \times 6 + 10 \times 10 = 133$) examples are given.
Obviously, no two arcs corresponding to  $x$
may have a common tail, but may have common heads.
Furthermore, for a monotone non-decreasing  $x$
the corresponding arcs cannot cross.
The same is true for  $y$.
However, two arcs corresponding to $x$  and to  $y$ may cross;
also the head of one may coincide with the tail of the other.

\begin{remark}
In contrast, for a monotone non-increasing mapping,
every two of its arcs intersect: either cross or have common heads.
Indeed, one turns a monotone non-decreasing mapping  $x$
(respectively,  $y$)  into a monotone non-increasing one  $x'$
(respectively,  $y'$) just replacing the order  of  $B$
(respectively,  of  $A$)  by the inverse one.
\end{remark}

By construction, digraph  $\Gamma$  is bipartite, with parts  $A$  and $B$.
Hence, every directed cycle (dicycle) in  $\Gamma$  is even.
There is an obvious one-to-one correspondence between
the dicycles of length $2$ (or  $2$-{\em dicycles}, for short) in $\Gamma(x,y)$
and the deals of  $G(x,y)$.
In the figures all $2$-{\em dicycles} are replaced by the nondirected edges.

\begin{proposition}
\label{p0}
Digraph  $\Gamma(x,y)$  contains at least one dicycle of length  $2$
and no longer dicycles.
\end{proposition}

\proof.
For any initial vertex  $v \in A \cup B$,
strategies  $x$  and  $y$  uniquely
define an infinite walk from  $v$.
Since sets $A$  and  $B$ are finite, this walk is a {\em lasso},
that is, it consists of an initial directed path and
a dicycle  $C$  repeated infinitely.
Furthermore,  $C$  is a $2$-dicycle whenever mappings
$x$  and  $y$  are monotone non-decreasing.
Indeed, it is easily seen that
either  $x$,  or  $y$, or both are not monotone whenever  $C$  is longer, since
in this case crossing arcs corresponding either to  $x$  or  to  $y$  must appear.
\qed


\medskip

The next statement results immediately from Proposition \ref{p0}.

\begin{proposition}
\label{p1}
Game correspondence $G$ is well-defined, that is,
$G(x,y) \neq \emptyset$  for all monotone non-decreasing
$x \in X$  and  $y \in Y$.
\qed
\end{proposition}

As Figures 1,\, 2, and 3 
show,
$G(x,y)$  may contain several deals.
These deals can be naturally ordered
in accordance with their importance, as the following claim shows.

\begin{proposition}
\label{p1a}
For any two deals  $(a,b), (a',b') \in G(x,y)$  we have:
either  $a > a'$  and  $b > b'$,  or  $a < a'$  and  $b < b'$.
\end{proposition}

\proof.
Indeed, equalities  $a=a'$ or $b=b'$  cannot hold,
because any two arcs corresponding to a mapping
cannot have a common tail.
Also they cannot cross whenever the considered mapping is monotone non-decreasing.
The above two observations imply the statement.
\qed

\begin{proposition}
\label{p1b}
Any set of pairwise disjoint
(nondirected)  edges between  $A$  and $B$  may form  $G(x,y)$
for some  monotone non-decreasing  $x$  and  $y$.
\end{proposition}

\proof.
Consider any  $k$  deals
$(a_{i_1}, b_{j_1}, \ldots, a_{i_k}, b_{j_k}$  such that
$a_{i_1} <  \ldots < a_{i_k}$  and   $b_{j_1} <  \ldots < b_{j_k}$. Define
$x(a) = b_{j_t}$  for all  $a$  such that  $a_{i_{t-1}} < a \leq  a_{i_t}$  and
$x(b) = a_{i_t}$  for all  $b$  such that  $b_{j_{t-1}} < b \leq  b_{j_t}$, where
the first inequalities are redundant for  $t=1$.
By construction,  $x$  and  $y$  are monotone non-decreasing and
$G(x,y) = \{(a_{i_1}, b_{j_1}, \ldots, a_{i_k}, b_{j_k}\}$.
\qed

\medskip

The remaining claims of this subsection we formulate and prove for Alice.
Of course, similar statements and proofs  hold for Bob as well: just replace
$x, X, a, A, i, I,$  "row", "Alice", and "she"  with
$y, Y, b, B, j, J,$  "column", "Bob", and "he", respectively, and vice versa.

\begin{proposition}
\label{p1c}
Any set of pairwise disjoint (nondirected) edges
between  $A$  and  $B$
corresponding to a monotone non-decreasing  $x$
forms   $G(x,y)$  for some  $y$.
\end{proposition}

\proof.
Define (a monotone non-decreasing) mapping $y$  as in the proof of Proposition \ref{p1b}.
\qed

\medskip

Choosing a game form  $g \in G$  we select one deal
for each pair $(x, y)$.
The following statement implies that, in a sense,
this choice does not matter.

\begin{proposition}
\label{p2}
For any deal  $(a,b) \in G(x, y)$ there exists
a strategy  $x' \in X$  such that  $G(x', y)$
consists of a unique deal  $\{(a,b)\}$.
\end{proposition}

\proof.
Obviously, the constant mappings
$x': A \rightarrow \{b\}$
satisfies the condition and proves the statement.
\qed

\medskip


Given a subset of outcomes  $W \subseteq O$  and
a game correspondence  $G : X \times Y \rightarrow 2^O$
(respectively, a game form  $g : X \times Y \rightarrow O$)
we say that Alice is effective  for  $W$  if she has
a strategy  $x \in X$  such that  $G(x, y) \subseteq W$
(respectively,  $g(x, y) \in W$)  for any  $y \in Y$;
or in other words, if there is a row  $x$
in which all outcomes are from  $W$.
We write $E(A, W) = 1$  in this case and  $E(A, W) = 0$  otherwise.

For each $W \subseteq O$, the following simple greedy algorithm verifies
in time linear in $|O| = mn$  whether Alice is effective for  $W$,
that is, if  $E(A, W) = 1$.
For each  $a_i \in A, \; i \in I = \{1, \ldots, m\}$,
select the minimal
$b_{j_i} \in B, \; j_i \in J = \{1, \ldots, n\}$,  if any,  such that
$(a_i, b_{j_i}) \in W$  and  $b_{j_i}  \geq  b_{j_{i-1}}$.
(The last condition is redundant for  $i=1$.)
If  such  $b_{j_i}$  exists for each  $i \in I$ then
stop and output  $E(B, W) = 1$.
In this case   $x : A \rightarrow B$  such that
$x(a_i) = b_{j_i}$  for all  $i \in I$  as a required strategy;
furthermore, by construction,  $x$  is the lexicographically minimal
monotone non-decreasing one. 
If  $b_{j_i}$  fails to exist for some  $i \in I$
then stop and output $E(A, W) = 0$.

\medskip

However, the last option is impossible when
each column contains an outcome from  $W$.

\begin{theorem}
\label{t1}
Given  $G : X \times Y \rightarrow 2^O$ and  $W \subseteq O$,
if for each  $y \in Y$  there exists an  $x \in X$  such that
$W \cap G(x,y) \neq \emptyset$  then
there exists an $x \in X$  such that
$G(x,y) \subseteq W$  for each  $y \in Y$.
\end{theorem}

This claim can be equivalently reformulated in many ways;
for example,  $E(B, O \setminus W) = 0 \; \Rightarrow \; E(A, W) = 1$;
see Section \ref{s2} for more options.

\medskip
 
\proof $\;$ of the theorem.
We apply the greedy algorithm verifying that
for each  $i = 1, \ldots, m$  a required
$b_{j_i}$  exists for  $a_i$.

First, let us consider the following simple case, just as an example.
Suppose that   $(a_1, b_2), (a_2, b_1) \in W$,
while  $(a_1, b_1) \not\in W$ and
$(a_2, b_j) \not\in W$  for all  $j = 2, \ldots, n$.
Then, $E(A, W) = 0$, since  $b_{i_2}$  does not exist.
However, in this case  $W$  cannot satisfy the conditions of the theorem.
Indeed, let us consider a
(monotone non-decreasing) strategy  $y \in Y$  such that
$y(b_1) = a_1$  and  $y(b_j) = a_2$  for all $j = 2, \ldots, n$.
Then, either  $(a_1, b_1) \in W$  or  $(a_2, b_j) \in W$
for at least one  $j = 2, \ldots, n$.
Thus, $G(x, y) \cap W = \emptyset$  for all  $x \in X$,
which contradicts the requirements to  $W$.

The general case can be analyzed similarly.
Suppose that the greedy algorithm assigns
$b_{j_i}$  to  $a_i$  for  $i = 1, \ldots, k$, for some  $k < m$,
but for  $a_{k + 1}$  there exists no
$b = b_{j_{k + 1}}$  satisfying the required conditions:
(l)  $(a_{k+1}, b) \in W$   and  (ll) $b \geq  b_{j_k}$,
that is, either (l), or (ll), or both fail for every $b$.
In each case we will construct a strategy  $y \in Y$
that violate conditions of the theorem.
When (l) fails, it is enough to take the constant mapping  $y$
that assigns  $b$  to every $a \in A$.
When (ll) fails,  $y$ is chosen in a more sophisticated way.
Recursively,  for  $j = 1, \ldots, k$, assign to each $b_j \in B$
the minimal  $a_{i_j} \in A$  such that
$(a_{i_j}, b_j) \not\in W$  and  $a_{i_j} > a_{i_{j - 1}}$.
In particular,   $a_{i_j} =  a_{k+1}$  for all  $j \geq j_k$.

By construction, $y$ is a monotone non-decreasing mapping
such that  $W \cap G(x,y) = \emptyset$  for all  $x \in X$,
which contradicts conditions of the theorem.
\qed

\section{Tightness and Nash-solvability}
\label{s2}

\subsection{Dual hypergraphs}
\label{ss20}
Given a ground-set  $O = \{o_1, \ldots, o_p\}$,
consider two hypergraphs
$$\cC = \{C_x \mid x \in X\} \; \text{and} \; \cD = \{D_y \mid 2^Y\}  \subseteq 2^O.$$

Here  $X, Y,$ and $O$  are arbitrary finite sets.
However, later they will be interpreted
as strategies  and outcomes, as in Sections \ref{s0} and \ref{s1}.
For this reason we keep the notation.
Consider the following properties  of  $\cC$  and  $\cD$:

\begin{itemize}
\item{($i$)} $C_x  \cap D_y \neq \emptyset$  for all  $x \in X$ and  $y \in Y$;
\item{($t'$)}  For any  $D' \subseteq O$  such that  $D' \cap C_x \neq \emptyset$
for each  $x \in X$  there exists an  $y \in Y$  such that  $D_y \subseteq D'$;
\item{($t''$)}  For any  $C' \subseteq O$  such that  $C' \cap D_y \neq \emptyset$
for each  $y \in Y$  there exists an  $x \in X$  such that  $C_x \subseteq C'$;
\item{($t$)} $C' \cap D' \neq \emptyset$  whenever
$C' \cap D_y \neq \emptyset$  for each $y \in Y$  and
$D' \cap C_x \neq \emptyset$  for each $x \in X$;
\item{$(r)$} For any  $x \in X, y \in Y$  and $o \in C_x \cap D_y$
there exist  an  $x' \in X$  and  $y' \in Y$  such that
$C_{x'} \cap D_y = C_x \cap D_{y'} = \{o\}$.
\end{itemize}

\begin{theorem}
\label{t2}
The following implications hold:
\begin{equation}
\label{eq2}
((i) \; \& \; (t))  \Leftrightarrow  ((i) \; \& \; (t'))
\Leftrightarrow ((i) \; \& \; (t'')) \Rightarrow  \; (r).
\end{equation}
\end{theorem}

Hypergraphs  $\cC$  and  $\cD$  satisfying the first three
equivalent properties are called {\em dual}.
Note that  $((i) \; \& \; (r))$  does not imply duality yet.

The last implication of Theorem \ref{t2} is obvious and all are well-known;
see, for example, \cite[Part I Section 4.2]{CH11},
where the above properties are formulated in the Boolean language,
in terms of the monotone
(positive) Conjunctive and Disjunctive Normal Forms;
see also \cite{Gur75}, \cite{Gur88}.
We leave to the reader to verify all implication of Theorem \ref{t2}.

\smallskip

It is also well-known \cite{CH11}  that
duality  of  $\cC, \cD \subseteq 2^O$
is respected by any identification
of the elements of  $O$.
Let us fix a subset  $O' \subseteq O$
such that  $|O'| > 1$  and
replace all  $o  \in O'$  by one new element $o'$.
This operation turns  $\cC$  and  $\cD$
into  $\cC'$  and  $\cD'$ , which remain dual
whenever  $\cC$   and  $\cD$  were dual.
Clearly, such operation can be applied successively several times.

\begin{remark}
A hypergraph is called {\em Sperner} if its
edges do not contain one another and,
in particular, cannot be equal.
Sperner hypergraphs are associated with
the  {\em irredundant} monotone
Conjunctive and Disjunctive Normal Forms \cite{CH11}.
For any hypergraph  $\cC$
(respectively, $\cD$) properties  $(i)$   and  $(t')$
(respectively,  $(i)$   and  $t'')$)  uniquely determine
the so-called  {\em dual} hypergraph  $\cC^d$
(respectively, $\cD^d$),  which is Sperner, by the definition.
If  $\cC$  (respectively, $\cD$) is Sperner itself then
duality is an involution, that is, $\cC^{dd} = \cC$
(respectively, $\cD^{dd} = \cD$).

Although in all our examples 
only dual pairs of Sperner hypergraphs appear,
but we do not assume that this always holds;
in particular, we allow hypergraphs to have embedded or equal edges.
\end{remark}

\subsection{Tight game correspondences and game forms}
\label{ss21}
Now, we interpret  $X, Y$ and $O$
as the sets strategies and outcomes.
A mapping  $G : X \times Y \rightarrow 2^O$
assigning a set of outcomes to each pair of strategies
is called a {\em game correspondence}.
A mapping  $g : X \times Y \rightarrow O$
assigning a single outcome of $O$  to each pair  $(x,y)$ 
is called a {\em game form}.
We will call  $G' : X \times Y \rightarrow 2^O$ a subcorrespondence
of  $G$  and write  $G' \leq G$  if  $G'(x,y) \subseteq G(x,y)$
for each  $x \in X$  and  $y  \in Y$.
If  $G' = g$  is a game form, we call it
a subform of  $G$  and write  $g \in G$.
All these concepts were already considered in Sections \ref{s0} and \ref{s1}
in the special case of monotone  bargaining;
here we extend them to the general case.

To each game correspondence  $G : X \times Y \rightarrow 2^O$
we assign two hypergraphs

$$
\cC = \{C_x = \bigcup_{y \in Y} G(x,y) \mid x \in X\}, \; 
\cD = \{D_y = \bigcup_{x \in X} G(x,y) \mid y \in Y\}  \subseteq 2^O  and 
$$
and call $G$ {\em tight} if they are dual.
Note that
$G(x,y) \subseteq  C_x \cap  D_y$  for every  $x \in X, y  \in Y$,
by definition, and hence,
property ($i$) holds
for   $\cC$  and  $\cD$ automatically, since  $G(x,y) \neq \emptyset$.

\smallskip

Game correspondences  $G_1, G_2, G_3, G_4$  in Figures 1--4
have the following four pairs of hypergraphs:

\smallskip

$\cC_1 = \{(o_{11}, o_{21}), (o_{11}, o_{22}),(o_{12}, o_{22})\},$

\smallskip

$\cD_1 = \{(o_{11}, o_{12}), (o_{11}, o_{22}), (o_{21}, o_{22})\};$

\medskip

$\cC_2 = \{(o_{11}, o_{21}, o_{31}), (o_{11}, o_{21}, o_{32}),
          (o_{11}, o_{22}, o_{32}), (o_{12}, o_{22}, o_{32})\},$

\smallskip

$\cD_2 = \{(o_{11}, o_{12}), (o_{11}, o_{22}), (o_{11}, o_{32}),
          (o_{21}, o_{22}), (o_{21}, o_{32}), (o_{31}, o_{32})\};$

\medskip

$\cC_3 = \{(o_{11}, o_{21}, o_{31}), (o_{11}, o_{21}, o_{32}),
          (o_{11}, o_{21}, o_{33}), (o_{11}, o_{22}, o_{32}),
          (o_{11}, o_{22}, o_{33}),
          \newline
          \;\;\;\;\;\;\;\;\;\;\;  (o_{11}, o_{23}, o_{33}),
          (o_{12}, o_{22}, o_{32}), (o_{12}, o_{22}, o_{33}),
          (o_{12}, o_{23}, o_{33}), (o_{13}, o_{23}, o_{33})\},$

\smallskip

$\cD_3 = \{(o_{11}, o_{12}, o_{13}), (o_{11}, o_{12}, o_{23}),
          (o_{11}, o_{12}, o_{33}), (o_{11}, o_{22}, o_{23}),
          (o_{11}, o_{22}, o_{33}),
          \newline
          \;\;\;\;\;\;\;\;\;\;\; (o_{11}, o_{32}, o_{33}),
          (o_{21}, o_{22}, o_{23}), (o_{21}, o_{22}, o_{33}),
          (o_{21}, o_{32}, o_{33}), (o_{31}, o_{32}, o_{33})\};$

\medskip

$\cC_5 = \cD_5 = \{(o_0,o_1,o_6), (o_0,o_2,o_5) (o_0,o_3,o_4), (o_1,o_2,o_4), (o_1,o_3,o_5), (o_2,o_3,o_6), (o_4,o_5,o_6)\}$; 

\medskip

$\cC_6 = \cD_6 = W_k = \{(o_0, o_1), (o_0, o_2), \dots, (o_0, o_k), (o_1, o_2, \ldots, o_k)\}$  for  $k = 3$; 
 
\medskip

$\cC_7 = \{(o_1,o_2,o_3),(o_1,o_2,o_4),(o_1,o_3,o_4),(o_2,o_3,o_4)\}$; 

\smallskip 

$\cD_7 = \{(o_1,o_2), (o_1,o_3), (o_1,o_4), (o_2,o_3), (o_2,o_4), (o_3,o_4)\}$;

\medskip

$\cC_8 = \{(o_1, o_2), (o_3, o_4), (o_1, o_4, o_5), (o_2, o_3, o_5)\},$

\smallskip

$\cD_8 = \{(o_1, o_3), (o_2, o_4), (o_1, o_4, o_5), (o_2, o_3, o_5)\}.$

\bigskip

We leave to the careful reader to

\begin{itemize}
\item{(j)} verify that each pair contains two dual
(Sperner) hypergraphs, or in other words, that the eight 
game correspondences in Figures 1--8 are tight;
\item{(jj)} reformulate all statements of Section \ref{ss20}
related to tightness for game correspondences and forms;
\item{(jjj)} check that some of them were proven in Section \ref{s1}
for the case of monotone bargaining.
\end{itemize}

\begin{proposition}
\label{p5}
If  $G' \leq G$  then
$G'$  is tight whenever  $G$   is tight;
moreover, in this  case $G$  and  $G'$
have the same pair of dual hypergraphs.
\end{proposition}

Note that the inverse implication fails.
For example, $G(x,y) = O$  may hold
for all  $x \in X$  and  $y \in Y$, in which case
$G$  is not tight unless
$|X| = 1$  or  $|Y| = 1$ or  $|O| = 1$, while
$G'$  may be tight.
Note also that $G'$  may be a game form.

\proof $\;$ of Proposition \ref{p5}.
If  $G$  is tight then property  $(r)$  holds.
implying that  $G$  and  $G'$  generate the same
pairs of hypergraphs.
\qed

\smallskip

Any two hypergraphs  $\cC$  and   $\cD$  on the ground-set  $O$
uniquely define the game correspondence
$G = G_{\cC, \cD}$  by formula  $G(x,y) = C_x \cap  D_y$.
By definition,  $G$  is tight if and only if
$\cC$   and   $\cD$  are dual.
By Proposition \ref{p5}, in this case every subcorrespondence
$G' \leq G$  and, in particular, each game form
$g \in G$  is tight too;
also tightness is preserved by any identifications of the outcomes,
as we already mentioned.
Note finally that hypergraphs  $\cC$  and  $\cD$  need not to be Sperner:
each of them may contain embedded and multiple edges.

\subsection{Nash-solvability}
\label{ss23}
To a game form  $g : X \times Y \rightarrow O$,
we assign a (normal form) game  $(g,u)$  introducing
the utility function  $u = (u_A, u_B)$  for the players:
$u_A : O \rightarrow \R$  for Alice and
$u_B : O \rightarrow \R$  for Bob. A strategy profile
$(x,y) \in X \times Y$  of a game  $(g,u)$  is called
a {\em Nash equilibrium} (NE)  if

\medskip
\noindent
$u_A(g(x',y)  \leq  u_A(g(x,y)$  for all $x' \in X$;
$\;\; u_B(g(x,y')  \leq  u_A(g(x,y)$  for all $y' \in Y$.

\medskip

A game   $(g,u)$  and its utility function  $u = (u_A, u_B)$
are called {\em zero-sum} if
$u_A(o)  + u_B(o) = 0$  for each  $o \in O$.
In this case a Nash equilibrium is usually referred to
as a {{\em saddle point}. A game form  $g$  is called
\begin{itemize}
\item{Nash-solvable (NS)}
if game  $(g,u)$  has a Nash equilibrium (NE) for any $u$;
\item{zero-sum-solvable}
if $(g,u)$  has a saddle point for every zero-sum  $u$;
\item{$\pm 1$-solvable}
if  $(g,u)$  has a saddle point for every
zero-sum  $u$  taking values $+1$  and  $-1$  only.
\end{itemize}

\begin{theorem}
\label{t-t}
(\cite{Gur75}, \cite{Gur88}).
The  following properties of game forms ar equivalent:
(S1) NS, (S2) zero-sum-solvability, (S3) $\pm 1$-solvability, (T) tightness.
\end{theorem}

Implications $(S1) \Rightarrow (S2) \Rightarrow (S3)$  are immediate from the definitions.
\newline
For the zero-sum case, the equivalence
$(T) \Leftrightarrow (S2) \Leftrightarrow (S3)$  was proven
by Edmonds and Fulkerson
(who called it the ``Bottle-Neck Extrema Theorem") in \cite{EF70}  and later,
independently, in \cite{Gur73}.  Then, NS
(and the implication $(T) \Rightarrow (S1)$)
was added to the list in \cite{Gur75}.
The same proof was repeated in a more focused paper \cite{Gur88}.
This proof is constructive:
a NE in the game $(g,u)$ is obtained by an explicit algorithm
whenever the game form  $g$  is tight.
Unfortunately, this algorithm is exponential in the size of $O$.

Another proof, based on a polynomial algorithm, appeared later;
see  \cite[Section 4.3]{Gur97}  and  \cite[Section 3]{BG03}.
Some new ideas from \cite{BBM89} and  \cite{DS91} were used in this proof.
Here, for reader's convenience, we reproduce it
in a stronger and slightly simplified form.

\proof $\;$ of the implication  $(T) \Rightarrow (S1)$.
We  will prove a stronger claim.
Given a tight game correspondence
$G : X \times Y \rightarrow  2^O$
and an arbitrary utility function  $u = (u_A, u_B)$,
we suggest an algorithm that finds
a simple NE
($(x,y)$  such that  $|G(x,y)| = 1$) in time linear in $|O|$.

Let us recall notation.
Given a subset of outcomes $W \subseteq O$,
we say that Alice is effective  for  $W$  if she has
a strategy  $x \in X$  such that  $G(x, y) \subseteq W$ for any  $y \in Y$;
or in other words, if there is a row  $x$
in which all outcomes are from  $W$, that is,  if  $C_x \subseteq W$.
We write $E(A, W) = 1$  in this case and  $E(A, W) = 0$  otherwise.
The effectivity function  $E(B, W)$  of Bob is  introduced similarly.

Obviously, the equalities
$E(A, W) = E(B, O \setminus W) = 1$
for no  $W \subseteq O$  can hold simultaneously.
Indeed, if Alice is effective for  $W$
and Bob for  $O \setminus W$  then
$\emptyset = G(x,y) \subseteq  C_x \cap D_y$
for the corresponding strategies  $x \in X$  and  $y \in Y$,
but  $G(x,y)$  cannot be empty.
Thus, implication
$E(A,W) = 1 \Rightarrow E(B, O \setminus W) = 0$
holds for any  $G$  and, by definition,
$G$  is tight if and only if the inverse implication holds:
$E(A,W) = 0 \Rightarrow E(B, O \setminus W) = 1$.

Consider partition  $O = W \cup W_A \cup W_B$
satisfying the following properties:

\begin{itemize}
\item{(a1)}
any outcome of  $W_A$  is not better for Alice
than any outcome of  $W$, that is,
$u_A(o') \leq u_A(o)$  for any  $o \in W$ and $o' \in W_A$;
\item{(b1)}
any outcome of  $W_B$  is not better for Bob
than any outcome of  $W$, that is,
$u_B(o') \leq u_B(o)$  for any  $o \in W$ and $o' \in W_B$;
\item{(a2)} $E(A, W_B) = 0$, that is, Alice cannot punish Bob;
\item{(b2)} $E(B, W_A) = 0$, that is, Bob cannot punish Allice.
\end{itemize}

Our arguments will be based on a dynamical process.
We begin with  $W = O$  and then reduce  $W$
sending outcomes to  $W_A$  and to $W_B$.
First, let us show that  $W$  cannot be eliminated completely.
Indeed, suppose that
$W = \emptyset$  and, hence,  $O = W_A \cup W_B$ is a partition.
Then, (a2), (b2), and tightness of $G$ imply that
$E(A, W_B) = E(B, W_A) = 1$, which is impossible,
as we already mentioned.
Thus, the dynamical process of reducing  $W$
will get stuck at some point, say, at a partition
$O = W \cup W_A \cup W_B$.

Let  $o^* \in O$  be the worst outcome in $W$  for Alice, that is,
$u_A(o^*) \leq u_A(o)$  for each $o \in W$.
We would like to move  $o^*$  from  $W$  to  $W_A$  but cannot.
The only possible reason is that (b2) will be violated,
that is, $E(B, W_A \cup \{o^*\}) = 1$.
Denote by  $y^*$  a strategy of Bob that
is effective for $(W_A \cup \{o^*\})$.

Moving  $o^*$  alone from to  $W$ to $W_B$  may violate (b1).
To keep it, we have to add to  $o^*$  all
outcomes  $W_B(o^*) \subseteq W$  that are not better than
$o^*$  for Bob,  that is,
$u_B(o^*) \geq u_B(o)$  for all $o \in W_B(o^*)$;
in particular, $o^* \in  W_B(o^*)$.
Then, moving $W_B(o^*)$  from  $W$  to $W_B$  is OK with (b1),
but still cannot be done.
The only possible reason is that (a2) is violated,
that is, $E(A, W_B \cup W_B(o^*)) = 1$.
Denote by  $x^*$  a strategy of Alice
that is effective for  $(W_B \cup W_B(o^*))$.

By construction,
$G(x^*, y^*) = (W_A \cup \{o^*\}) \cap (W_B \cup W_B(o^*)) = \{o^*\}$.
and  $(x^*, y^*)$   is a NE, moreover, it is simple,
since $G(x^*, y^*)$  consists of the single outcome $o^*$.
\qed

\section{More classes of tight game forms}
\label{s3}
A criterion of tightness for the game forms
associated with veto voting schemes can be found in the textbooks
\cite[Chapter 6]{Mou83} or \cite[Chapter 6]{Pel84}; see also \cite{GG86, Gur17a}.
It can be reformulated as a criterion of duality for the monotone
(positive) threshold Boolean functions \cite[Part II Chapter 9]{CH11}.

Here we will recall five more classes of tight game forms, whose
NS follows in all cases from Theorem \ref{t-t}, but sometimes
can be shown simpler.

\subsection{Positional game structures}
\label{ss30}
Let $\Gamma = (V,E)$ be a digraph whose vertices
and arcs are interpreted as positions and moves, respectively.
Furthermore, let  $V_A$  and  $V_B$  be positions
of positive out-degree, controlled by Alice and Bob, respectively,
while $V_T$  be the set of terminal positions, of out-degree zero.
We assume that  $V = V_A \cup  V_B \cup V_T$  is a partition.
A strategy $x \in X$  of Alice
(resp., $y \in Y$  of Bob) is a mapping that
assigns to each position  $v \in V_A$
(resp., $v \in V_B$) an arbitrary move from this position.
An initial position  $v_0 \in V_A \cup V_B $  and a pair
$(x,y)$  of strategies uniquely define a walk
that begins and  $v_0$  and then follows
the decisions made by  $x$  and $y$.
In theory of positional games this walk is called a {\em play}.
Each play either terminates in  $V_T$  or it is infinite;
in the latter case, it must form a lasso:
first, an initial path and then a dicycle repeated infinitely.

Positional structures defined above
can also be represent in normal form.
We will introduce a game form  $g : X \times Y \rightarrow O$,
where standardly  $O$  denotes a set of outcomes.
There are several ways to introduce this set.
One is to identify all infinite
plays and consider them as a single outcome  $c$;
then  $O = V_T \cup \{c\}$.
Such games were introduced by Washburn \cite{Was90}
who called them {\em deterministic graphical} (DG) games.
In this case, for any positional structure
the associated game form is tight
(\cite[Section 3]{BG03}; see also \cite[Section 12]{BGMS07})
and, hence, it is also NS, as we know).

The following generalization was suggested recently in \cite{Gur17}.
Digraph  $\Gamma$  is called {\em strongly connected}  if
for any  $v, v' \in V$  there is a directed path from  $v$  to  $v'$
(and, hence, from  $v'$ to $v$, as well).
By this definition, the union of two strongly connected digraphs
is strongly connected  whenever they have a common vertex.
A vertex-inclusion-maximal strongly connected induced subgraph of $\Gamma$
is called its {\em strongly connected component}.
Obviously, any digraph  $\Gamma = (V, E)$  admits a unique decomposition into such components:
$\Gamma^j = \Gamma[V^j] = (V^j, E^j)$  for $j \in J$, where
$J$ is a set of indices; furthermore
$V = \cup_{j \in J} V^j$  is a partition of  $V$.
It can be determined in time linear in
the size of  $\Gamma$ (that is, in $(|V| + |E|)$) and
has numerous applications; see \cite{Tar72} and \cite{Sha81} for more details.

One more application, to positional games was suggested in \cite{Gur17}.
For each  $j \in J$, let us contract  $\Gamma^j$  into a single vertex  $v^j$.
Then, all edges of  $E^j$
(in particular, the loops) disappear and
we obtain an acyclic digraph  $\Gamma^* = (V^*, E^*)$.

In the DG games, all infinite plays (lassos) represent a unique outcome  $c$.
In \cite{Gur17}, the following weaker assumption was considered.
Two plays are equivalent, and represent the same outcome,
when the directed cycles of their lassos are contained
in the same strongly connected  component of  $\Gamma$.
Thus,  $J$  becomes the set of outcomes.
Such games were called {\em deterministic graphical multi-stage} (DGMS).
In this case too, for any positional structure,
the associated game form is tight \cite{Gur17}.
This result DGMS games strengthens
the similar one for the DG games, because
merging outcomes, obviously, respects tightness and NS.

Every strongly connected component  $\Gamma^j$
can be interpreted as a diplomatic issue.
(The level may be not very high; see examples in \cite{Gur17}).
If this issue is resolved, the play leaves $\Gamma^j$,
enters another component,
and the parties (players) proceed with the negotiation process.
Otherwise, the play cycles in  $\Gamma^j$  and negotiations terminate
resulting in the outcome  $c_j$.
It only matters at what issue  $\Gamma^j$
the parties have been forced into an impasse, while details
(a particular lasso, or a directed cycle in which it ends in $\Gamma^j$)  are irrelevant.
Thus, a DGMS game is interpreted as a multistage diplomatic negotiations.
The main result shows that the corresponding DGMS  game form is NS;
furthermore, for any payoff a NE in the obtained game can be found in linear time,
in other words, any bilateral multistage diplomatic conflict can be efficiently resolved
and the resolving strategies can be efficiently determined.

\bigskip

One can also treat each dicycle of  $\Gamma$
as a separate outcome; then,  $O = V_T \cup C(\Gamma)$.
In this case, not all but only some special game forms are tight.
Their characterization is not known, in general,
but the problem is solved for the symmetric game structures, that is,
assuming that  $(v',v) \in E$   whenever  $(v,v') \in E$.
First, tight game forms of the symmetric bipartite game structures
were characterized in  \cite{GG91}, \cite{GG91a}, and \cite{GG92}.
Then, in \cite{BGMS07} this result was strengthened by waving
the assumption that  $\Gamma$  is bipartite.

\subsection{Selfdual hypergraphs}
\label{ss31}

A hypergraph  $\cC$  is called {\em selfdual}
if it is dual to itself.
Obviously, each such hypergraph generates
a symmetric tight game correspondence.
The Fano projective plane provides an example:
$$\cC = \{(o_0, o_1, o_6), (o_0, o_2, o_5), (o_0, o_3, o_4),
(o_1, o_2, o_4), (o_1, o_3, o_5), (o_2, o_3, o_6),
(o_4, o_5, o_6)\}.$$

\bigskip

\begin{tabular}{c||c|c|c|c|c|c|c}
& $o_0 o_1o_6$ & $o_0 o_2o_5 $ & $o_0 o_3 o_4$ & $o_1 o_2 o_4 $ & $o_1 o_3o_5$& $o_2o_3o_6$ & $o_4o_5o_6$
\\ \hhline{--------}
$o_0 o_ 1 o_6$ & $o_0o_1o_6$ & $o_0$ & $o_0$ & $o_1$ & $o_1$ & $o_6$ & $o_6$
\\ \hhline{--------}
$o_0 o_2 o_5$ & $o_0$ & $o_0o_2o_5$ & $o_0$ & $o_2$ & $o_5$ & $o_2$ & $ o_5$
\\ \hhline{--------}
$o_0 o_3 o_4$ & $o_0$ & $o_0$ & $o_0 o_3 o_4$ & $o_4$ &$o_3$ & $o_3$ & $o_4$
\\ \hhline{--------}
$o_1 o_2 o_4$ & $o_1$ & $o_2$ & $o_4$ & $o_1o_2 o_4 $ & $o_1$ & $o_2$ & $o_4$
\\ \hhline{--------}
$o_1 o_3 o_5$ & $o_1$ & $o_5$ & $o_3$ & $o_1$ & $o_1o_3 o_5 $ & $o_3$ & $o_5$
\\ \hhline{--------}
$o_2 o_3 o_6$ & $o_6$ & $o_2$ & $o_3$ & $o_2$ & $o_3$ & $o_2o_3 o_6 $ & $o_6$
\\ \hhline{--------}
$o_4 o_5 o_6$ & $o_6$ & $o_5$ & $o_4$ & $o_4$ & $o_4$ & $o_6$ & $o_1o_2 o_4 $
\end{tabular}

\bigskip

Figure 4. Game correspondence of the Fano plane. 

\bigskip

To every pair of hypergraphs  $\cC$  and  $\cD$
defined on the same ground-set  $O$,
assign the hypergraph  $\cE = \{c \cC, d \cD, (c,d)\}$,
where $c$  and  $d$  are two new objects, $c, d \not\in O$. 
To get  $\cE$, we add  $c$ (resp., $d$)
to each edge of  $\cC$ (resp., of $\cD$) and
also add one extra edge $(c, d)$. 
Seymour noticed that  $\cC$  and  $\cD$  are dual if and only if
$\cE$ is selfdual  \cite{Sey74}.

\bigskip

Let us define a large class of selfdual hypergraphs.
Given  $p = k+1 \geq 3$  and the ground-set $\{o_0, o_1, \ldots, o_k\}$,
the $k$-{\em wheel} is defined as the hypergraph
$$W_k = \{(o_0, o_1), (o_0, o_2), \dots, (o_0, o_k), (o_1, o_2, \ldots, o_k)\}.$$
It is easily seen that  $W_k$  is selfdual for any  $k \geq 2$.

Consider the following interpretation.
Two players are choosing among $k+1$ options
$\{o_0, o_1, \ldots , o_k\}$, where
$o_0$  is a special one, which means that
players did not come to an agreement.
Each player has $k+1$  strategies, either to suggest
an option, except  $o_0$, or to say ``any".
The game correspondence is defined as follows:

\begin{itemize}
\item
If they suggest  $o_i$  and  $o_j$
(recall that  $i, j \in \{1, \ldots, k\}$)
then  $o_0$ is chosen when  $i \neq j$  and
either $o = o_i = o_j$  or  $o_0$  is chosen when  $i = j$.
\item
If one suggest  $o_i, \; i = 1,\ldots, k,\;$  while the other ``any",
then  $o_i$  is chosen.
\item
If both say ``any" then each option, except $o_0$, can be chosen.
\end{itemize}

For  $k=3$  the game correspondence is given on the following table.

\bigskip

\begin{tabular}{c||c|c|c|c}
& $o_0 o_1$ & $o_0 o_2 $ & $o_0 o_3$ & $o_1 o_2 o_3 $
\\ \hhline{-----}
$o_0 o_1$ & $o_0o_1$ & $o_0$ & $o_0$ & $o_1$
\\ \hhline{-----}
$o_0 o_2$ & $o_0$ & $o_0o_2$ & $o_0$ & $o_2$
\\ \hhline{-----}
$o_0 o_3$ & $o_0$ & $o_0$ & $o_1 o_3$ & $o_3$
\\ \hhline{-----}
$o_1 o_2 o_3$ & $o_1$ & $o_2$ & $o_3$ & $o_1o_2 o_3 $
\end{tabular}

\bigskip 

Figure 5. Game correspondence of the wheel, $\cC = \cD = W_3$.

\bigskip

For each $k \geq 2$, by Theorem \ref{t-t},
such game correspondence is NS.
Moreover, we proved that a simple NE exists for every payoff.
This proof can be significantly simplified in the considered case.

Let player $1$ (resp., $2$) chooses his (her)
best option  $o_i$ (resp., $o_j$)
among $o_1, \ldots, o_k$, while
player $2$ (resp., $1$)  says ``any".
Then  $o_i$   (resp., $o_j$) in the last column
(resp., row) will be a simple NE unless
$o_0$  is better for the opponent.
We can assume that the latter holds for both players, since
in all other cases we  are done.
If  $i \neq j$ then the intersection
of the row  $i$  and column  $j$ form a simple NE 
with the outcome $o_0$.
Interestingly, the most difficult case is  $i = j$,
when the same option is the best for both.
The corresponding diagonal entry is always a NE
with the outcome either  $o = o_i = o_j$  or  $o_0$,
but this NE is not simple.
Yet, the case is still simple.
Indeed, since $o_0$
is the best option for  both players,
any simple strategy profile with the outcome  $o_0$
(and there exist  $k(k-1)$  of them) is a NE.

\subsection{Symmetric dual hypergraphs}
\label{ss32}
Let  $k$ and  $\ell$  be two positive integers and 
$\cC$  and  $\cD$  be symmetric hypergraphs
of dimension  $k$ and  $\ell$  
on the ground-set $O = \{o_1, \ldots, o_p\}$  such that 
$p = k + \ell - 1$. 
In other words,  $\cC $  (resp., $\cD$)  consists of all
$k$-subsets  (resp., $\ell$-subsets) of $O$.
It is well-known that  $\cC$  and  $\cD$  are dual;
see, for example, \cite[Part I Chapter 4]{CH11}.

A natural interpretation can be given in terms of veto voting.
There are $p$  candidates and  two voters with  veto power  $p - k$  and  $p - \ell$;
see \cite[Chapter 6]{Mou83}, \cite[Chapter 6]{Pel84}, or \cite{GG86}.
For  $k=3$ and $\ell = 2$ the game correspondence is given on the following table.

\medskip

\begin{tabular}{c||c|c|c|c|c|c}
& $12$ & $13 $ & $14$ & $23 $ & $24$ & $34$
\\ \hhline{-------}
$123$ & $12$ & $13$ & $1$ & $23$ & $2$ & $3$
\\ \hhline{-------}
$124$ & $12$ & $1$ & $14$ & $2$ & $24$ & $4$
\\ \hhline{-------}
$134$ & $1$ & $13$ & $14$ & $3$ & $4$ & $34$
\\ \hhline{-------}
$234$ & $2$ & $3$ & $4$ & $23 $ & $24$ & $34$
\end{tabular}

\medskip

Figure 6. Game correspondence of the dual symmetric hypergraphs 
of  dimensions  $k = 3$ and $\ell = 2$  on the ground-set of size $p = 4$  
represents veto voting with  $p = 4$  candidates and 
$n = 2$  voters of the veto powers  $p - k  = 1$  and  $p - \ell = 2$, respectively.         

\medskip

For any positive integer  $k$ and $\ell$, a simple NE exists;
see the proof of Theorem \ref{t-t}.
It can be simplified in the considered special case.
Without loss of generality, we  may assume that
voter  $1$  prefers the candidates from $O$ in the order 
$u_1(o_1) \leq \ldots \leq u_1(o_p)$  and, moreover, that
all inequalities here are strict.
He can veto $k$  his worst candidates and guarantee  $o_{k+1}$ or better.
If the game is zero-sum, voter $2$
just veto her $\ell$  worst candidates $O_2 = \{o_{k+2}, \ldots, o_p\}$.
Obviously, these two strategies form a simple NE,
electing the centrist $o_{k+1}$.

If the game is not zero-sum then 
among $\ell$ candidates of   $O_2$  may exist  $r$ 
that are better than $o_{k+1}$  for voter  $2$
(and for voter  $1$  as well.
Note that  $r = 0$  in the zero-sum case).
Let us choose among these $r$ candidates
one with the largest index, say, $o_m$; 
in other words, the best one for both voters. 
By construction, there are two sets  $C, D \subseteq O$  such that
$|C| = k$, $|D| = \ell$, $C \cap D = \{o_m\}$, and $o_m$
is better than any $o \in C$  for  voter $1$
and any  $o \in D$  for voter  $2$.
Thus, strategy profile  $(C, D)$ forms a simple NE.

\subsection{Jordan dual hypergraphs.
Choosing Battlefields in Wonderland}
\label{ss33}
Wonderland is a subset of the plane
homeomorphic to the closed disc.
Without loss of generality, one can choose a square  $Q$
with the sides $N, E, S, W$, as in Figure 7.
Let us partition $Q$  into areas
$O = \{o_1, \ldots, o_p\}$ each of which
is homeomorphic to the closed disc, too.
Any two distinct areas  $o_i, o_j \in O$  are either disjoint or
intersect in a set homeomorphic to a closed interval
that contains more than one point.
Equivalently, we can require that the borders
of the areas of  $O$  form a regular graph of degree $3$,
as in Figure 7.

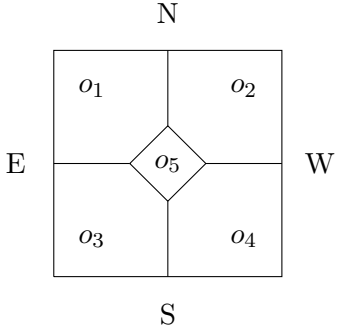
\begin{tikzpicture}
\draw (-1.5,-1.5) -- (-1.5,1.5) -- (1.5,1.5) -- (1.5,-1.5) -- (-1.5,-1.5);
\draw (1.5,0)--(0.5,0);\draw (-1.5,0)--(-0.5,0);\draw(0,1.5)--(0,0.5);\draw(0,-1.5)--(0,-0.5);
\draw(0.5,0)--(0,0.5)--(-0.5,0)--(0,-0.5)--(0.5,0);
\node at (0,0) {$o_5$};\node at (-1,-1) {$o_3$};\node at (1,-1) {$o_4$};\node at (-1,1) {$o_1$};\node at (1,1) {$o_2$};
\node at (2,0) {W};\node at (0,2) {N};\node at (-2,0) {E};\node at (0,-2) {S};
\end{tikzpicture}

\vskip 1cm

\begin{tabular}{c||c|c|c|c}
& $o_1 o_3$ & $o_2 o_4 $ & $o_1 o_4 o_5$ & $o_2 o_3 o_5 $
\\ \hhline{-----}
$o_1 o_2$ & $o_1$ & $o_2$ & $o_1$ & $o_2$
\\ \hhline{-----}
$o_3 o_4$ & $o_3$ & $o_4$ & $o_4$ & $o_3$
\\ \hhline{-----}
$o_1 o_4 o_5$ & $o_1$ & $o_4$ & $o_1 o_4 o_5$ & $o_5$
\\ \hhline{-----}
$o_2 o_3 o_5$ & $o_3$ & $o_2$ & $o_5$ & $o_2 o_3 o_5$
\end{tabular}

\bigskip 

Figure 7. The Jordan game correspondence of a map. 

\bigskip 

Two players, Tweedledee  and Tweedledum, agreed to have a battle.
The next thing to do is to choose a battlefield,
which should be an area  $o \in O$.
The strategies  $x \in X$  of Tweedledee
(resp., $y \in Y$  of Tweedledum) are all
inclusion-minimal subsets
$x \subseteq O$ (resp., $y \subseteq O$
connecting  $N$  and  $S$  (resp., $E$  and  $W$).
The game correspondence}  $G : X \times Y \rightarrow O$ 
defined by the equality  $G(x,y) = x \cap y \subseteq O$ 
will be called {\em Jordan};  
see again Figure 7 for an example.
The following ``topological" result is well-known.

\begin{theorem}
\label{tt-battle}
Any Jordan game correspondence $G$  is
well defined (that is, $G(x,y) \neq \emptyset$) and tight.
\end{theorem}

\proof.
It is easily seen that any two pathes in  $Q$,
one connecting  $N$  with  $S$  and
another connecting  $E$  with  $W$,  intersect.
This is one of the versions of the Jordan curve theorem.
(Note that the claim fails when 
vertices of degree larger than  $3$  are allowed in  $Q$.)  
Hence, $G(x, y) = x \cap y \neq \emptyset$  for any  $x$  and  $y$.
Moreover, it follows from the same theorem that
a subset of  $Q$  blocking all pathes between  $N$ and $S$
must contain a path between $E$ and $W$, and vice versa.
As we know, each of these two properties is
equivalent with tightness of  $G$.
\qed

\smallskip

Thus, all results of Sections \ref{s1} and \ref{s2} related to NS follow.
Let us note that in this case we do not know any simpler way than
to derive the existence of a simple NE from Theorem \ref{t-t}.

\begin{remark}
We required inclusion-minimality
for the subsets  $x, y \in O$  just to reduce the number of strategies.
This requirement can be waved, but 
Theorem \ref{tt-battle} with 
all its corollaries related to NS still hold.
\end{remark}

\section{Three-person monotone bargaining may be not tight and not NS}
\label{s4}
Let us introduce Claire to Alice and Bob
and consider three sets of objects
$$A = \{a_1, \ldots, a_m\},  B = \{b_1, \ldots, b_n\}, C = \{c_1, \ldots, c_k\},$$
\noindent
and the set of outcomes (or possible dealls)  $O = A \times B \times C$.

The strategies of Alice, Bob, and Clair are
all monotone non-decreasing mappings
$x : A \rightarrow B, y : B \rightarrow C$, and  $z : C \rightarrow A$.
Given a triplet of strategies  $(x,y,z) \in X \times Y \times Z$, naturally, a
triplet of objects  $(a,b,c) \in A \times B \times C$
is called a {\em deal} if  $x(a) = b, y(b) = c$, and  $z(c) = a)$.

\begin{proposition}
\label{p7}
Every triplet  $(x,y,z)$  generates at least one deal.
\end{proposition}

\proof.
We naturally assign to  $(x,y,z)$  a tripartite graph
$\Gamma = \Gamma(x,y,z)$  with the vertex-set  $V = A \cup B \cup C$
and all arcs   $(a,b), (b,c)$, and $(c,a)$  such that
$x(a) = b, y(b) = c$, and  $z(c) = a)$.
Then, any initial vertex  $v_0 \in V$
uniquely defines a walk.
Clearly, this walk is a lasso ending in a $3t$-dicycle;
moreover, $t=1$  follows from the fact that
$x,y$, and $z$  are monotone non-decreasing he mappings.
\qed

\smallskip

Thus, $G(x,y,z) \neq \emptyset$  for all triplets  $(x,y,z)$
and we obtain a game correspondence
$G : X \times Y \times Z \rightarrow 2^O$.
Yet, already for $m = n = k = 2$
the correspondence $G_{2,2,2}$ in Figure 8
is not tight and not  NS.
Indeed, it is not difficult to see that
Alice's and Bob-Claire's hypergraphs

\smallskip
\noindent
$\cH_A = \{
(o_{111}, o_{112}, o_{121}, o_{122}),
(o_{111}, o_{121}, o_{212}, o_{222}),
(o_{211}, o_{212}, o_{221}, o_{222})\};$

\noindent
$\cH_{BC} = \{
(o_{111}, o_{211}),
(o_{112}, o_{212}),
(o_{111}, o_{221}),
(o_{111}, o_{222}),
(o_{112}, o_{222}),
(o_{121}, o_{221}),
(o_{121}, o_{222})\}.$

\smallskip
\noindent
are not dual. For example,
$(o_{121}, o_{211})$  and  $(o_{121}, o_{212})$  are transversal to  $\cH_A$
that do not appear in $\cH_{BC}$.

\smallskip

NS fails too, for example, for any
utility  $u = (u_A, u_B, u_C)$  such that

\smallskip
\noindent
$u_A(212) > u_A(222) > u_A(111), \, u_A(121) > u_A(111), \, u_A(122) > u_A(112)$;

\smallskip
\noindent
$u_B(112) > u_B(111),  u_B(121) > u_B(122),  u_B(222) > u_B(221),   u_B(212) > u_B(211)$;

\smallskip
\noindent
$u_C(111) > u_C(222),  u_C(122) > u_C(122),  u_C(112) > u_C(212),   u_C(221) > u_C(121)$.

\bigskip

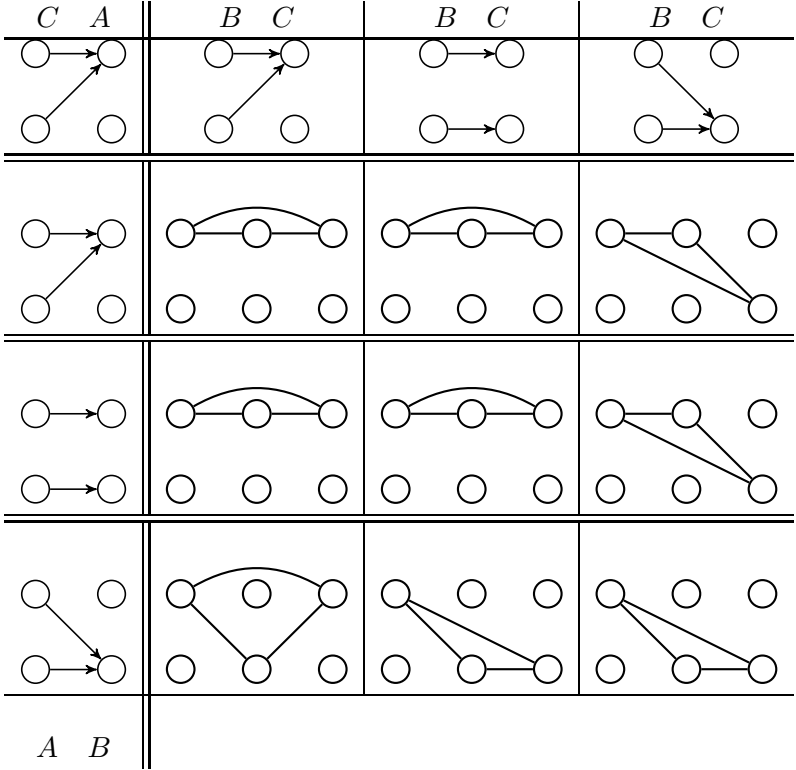
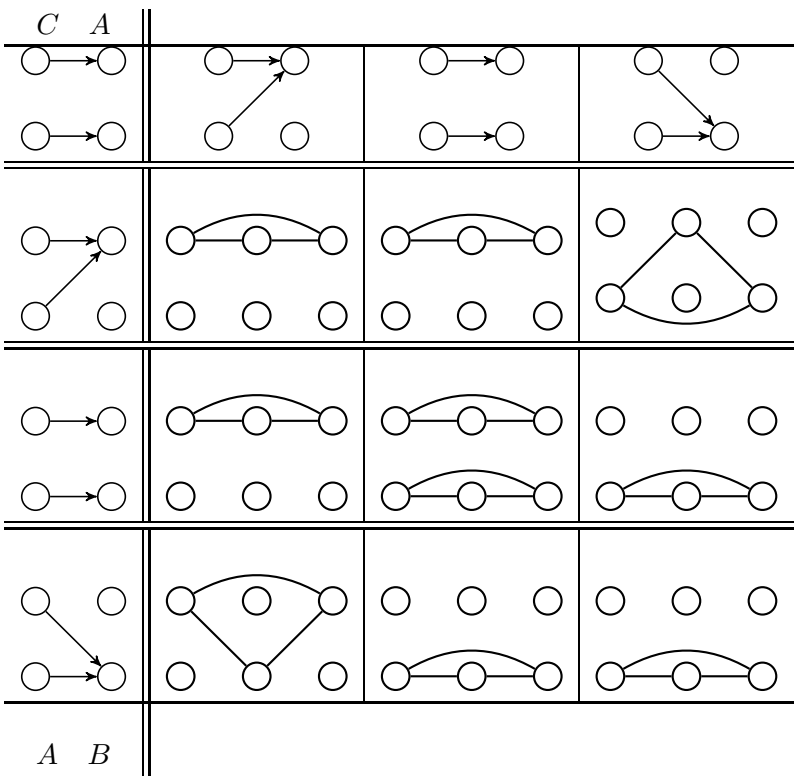
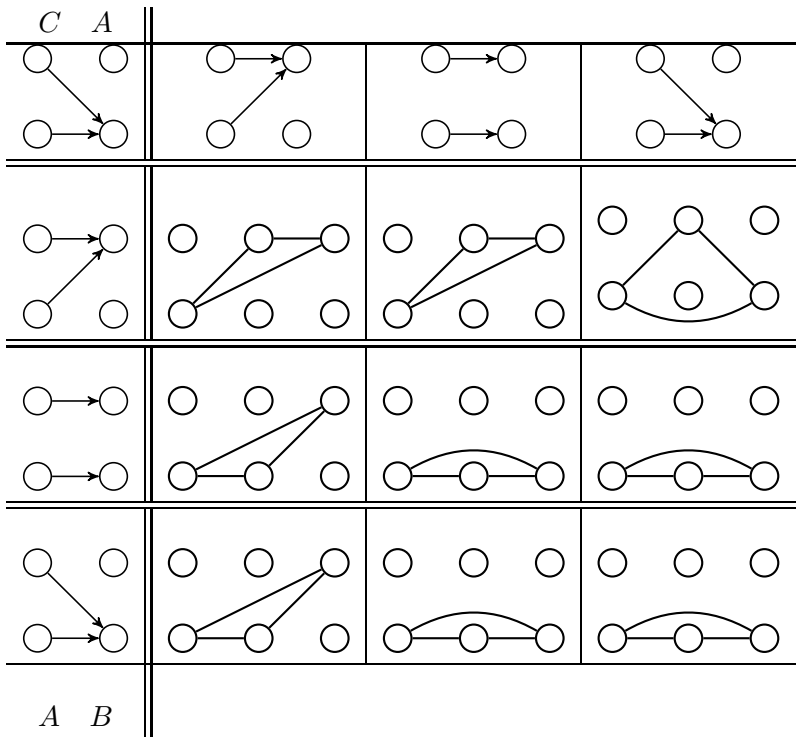
\begin{tabular}{c||c|c|c}
$C\quad A$&$B \quad C$&$B \quad C$&$B\quad C$
\\ \hhline{----}
\begin{tikzpicture}[>=stealth',semithick,auto,a/.style={circle,fill=white,draw},b/.style={circle,fill=white,draw}]\node[a] (00) at (0,0) {};\node[a] (10) at (0,-1) {};\node [b] (01) at (1,0) {};\node [b] (11) at (1,-1) {}; \path[->] (10) edge (01); \path[->] (00) edge (01); \end{tikzpicture} &
\begin{tikzpicture}[>=stealth',semithick,auto,a/.style={circle,fill=white,draw},b/.style={circle,fill=white,draw}]\node[a] (00) at (0,0) {};\node[a] (10) at (0,-1) {};\node [b] (01) at (1,0) {};\node [b] (11) at (1,-1) {}; \path[->] (00) edge (01); \path[->] (10) edge (01); \end{tikzpicture} &
\begin{tikzpicture}[>=stealth',semithick,auto,a/.style={circle,fill=white,draw},b/.style={circle,fill=white,draw}]\node[a] (00) at (0,0) {};\node[a] (10) at (0,-1) {};\node [b] (01) at (1,0) {};\node [b] (11) at (1,-1) {}; \path[->] (00) edge (01); \path[->] (10) edge (11); \end{tikzpicture} &
\begin{tikzpicture}[>=stealth',semithick,auto,a/.style={circle,fill=white,draw},b/.style={circle,fill=white,draw}]\node[a] (00) at (0,0) {};\node[a] (10) at (0,-1) {};\node [b] (01) at (1,0) {};\node [b] (11) at (1,-1) {}; \path[->] (00) edge (11); \path[->] (10) edge (11); \end{tikzpicture}
\\ \hhline{====} & & & \\[0.1em]
\begin{tikzpicture}[>=stealth',semithick,auto,a/.style={circle,fill=white,draw},b/.style={circle,fill=white,draw}]\node[a] (00) at (0,0) {};\node[a] (10) at (0,-1) {};\node [b] (01) at (1,0) {};\node [b] (11) at (1,-1) {}; \path[->] (00) edge (01); \path[->] (10) edge (01); \end{tikzpicture} &
\begin{tikzpicture}[thick,auto,a/.style={circle,fill=white,draw},b/.style={circle,fill=white,draw}]\node[a] (00) at (0,0) {};\node[a] (10) at (0,-1) {};\node [a] (01) at (1,0) {};\node [a] (11) at (1,-1) {};\node[b] (02) at (2,0){};\node[b] (12) at (2,-1){};\draw (00) edge (01);\draw (01) edge (02);\draw (00) edge[bend left] (02);\end{tikzpicture}  &
\begin{tikzpicture}[thick,auto,a/.style={circle,fill=white,draw},b/.style={circle,fill=white,draw}]\node[a] (00) at (0,0) {};\node[a] (10) at (0,-1) {};\node [a] (01) at (1,0) {};\node [a] (11) at (1,-1) {};\node[b] (02) at (2,0){};\node[b] (12) at (2,-1){};\draw (00) edge (01);\draw (01) edge (02);\draw (00) edge[bend left] (02);\end{tikzpicture}  &
\begin{tikzpicture}[thick,auto,a/.style={circle,fill=white,draw},b/.style={circle,fill=white,draw}]\node[a] (00) at (0,0) {};\node[a] (10) at (0,-1) {};\node [a] (01) at (1,0) {};\node [a] (11) at (1,-1) {};\node[b] (02) at (2,0){};\node[b] (12) at (2,-1){};\draw (00) edge (01);\draw (01) edge (12);\draw (00) edge (12);\end{tikzpicture}
\\ \hhline{====} & & & \\[0.1em]
\begin{tikzpicture}[>=stealth',semithick,auto,a/.style={circle,fill=white,draw},b/.style={circle,fill=white,draw}]\node[a] (00) at (0,0) {};\node[a] (10) at (0,-1) {};\node [b] (01) at (1,0) {};\node [b] (11) at (1,-1) {}; \path[->] (00) edge (01); \path[->] (10) edge (11); \end{tikzpicture} &
\begin{tikzpicture}[thick,auto,a/.style={circle,fill=white,draw},b/.style={circle,fill=white,draw}]\node[a] (00) at (0,0) {};\node[a] (10) at (0,-1) {};\node [a] (01) at (1,0) {};\node [a] (11) at (1,-1) {};\node[b] (02) at (2,0){};\node[b] (12) at (2,-1){};\draw (00) edge (01);\draw (01) edge (02);\draw (00) edge[bend left] (02);\end{tikzpicture}  &
\begin{tikzpicture}[thick,auto,a/.style={circle,fill=white,draw},b/.style={circle,fill=white,draw}]\node[a] (00) at (0,0) {};\node[a] (10) at (0,-1) {};\node [a] (01) at (1,0) {};\node [a] (11) at (1,-1) {};\node[b] (02) at (2,0){};\node[b] (12) at (2,-1){};\draw (00) edge (01);\draw (01) edge (02);\draw (00) edge[bend left] (02);\end{tikzpicture}  &
\begin{tikzpicture}[thick,auto,a/.style={circle,fill=white,draw},b/.style={circle,fill=white,draw}]\node[a] (00) at (0,0) {};\node[a] (10) at (0,-1) {};\node [a] (01) at (1,0) {};\node [a] (11) at (1,-1) {};\node[b] (02) at (2,0){};\node[b] (12) at (2,-1){};\draw (00) edge (01);\draw (01) edge (12);\draw (00) edge (12);\end{tikzpicture}
\\ \hhline{====} & & & \\[0.1em]
\begin{tikzpicture}[>=stealth',semithick,auto,a/.style={circle,fill=white,draw},b/.style={circle,fill=white,draw}]\node[a] (00) at (0,0) {};\node[a] (10) at (0,-1) {};\node [b] (01) at (1,0) {};\node [b] (11) at (1,-1) {}; \path[->] (00) edge (11); \path[->] (10) edge (11); \end{tikzpicture} &
\begin{tikzpicture}[thick,auto,a/.style={circle,fill=white,draw},b/.style={circle,fill=white,draw}]\node[a] (00) at (0,0) {};\node[a] (10) at (0,-1) {};\node [a] (01) at (1,0) {};\node [a] (11) at (1,-1) {};\node[b] (02) at (2,0){};\node[b] (12) at (2,-1){};\draw (00) edge (11);\draw (11) edge (02);\draw (00) edge[bend left] (02);\end{tikzpicture}  &
\begin{tikzpicture}[thick,auto,a/.style={circle,fill=white,draw},b/.style={circle,fill=white,draw}]\node[a] (00) at (0,0) {};\node[a] (10) at (0,-1) {};\node [a] (01) at (1,0) {};\node [a] (11) at (1,-1) {};\node[b] (02) at (2,0){};\node[b] (12) at (2,-1){};\draw (00) edge (11);\draw (11) edge (12);\draw (00) edge (12);\end{tikzpicture}  &
\begin{tikzpicture}[thick,auto,a/.style={circle,fill=white,draw},b/.style={circle,fill=white,draw}]\node[a] (00) at (0,0) {};\node[a] (10) at (0,-1) {};\node [a] (01) at (1,0) {};\node [a] (11) at (1,-1) {};\node[b] (02) at (2,0){};\node[b] (12) at (2,-1){};\draw (00) edge (11);\draw (11) edge (12);\draw (00) edge (12);\end{tikzpicture}
\\ \hhline{----} \\
$A\quad B$
\end{tabular}

\vskip1cm

\begin{tabular}{c||c|c|c}
$C\quad A$
\\ \hhline{----}
\begin{tikzpicture}[>=stealth',semithick,auto,a/.style={circle,fill=white,draw},b/.style={circle,fill=white,draw}]\node[a] (00) at (0,0) {};\node[a] (10) at (0,-1) {};\node [b] (01) at (1,0) {};\node [b] (11) at (1,-1) {}; \path[->] (10) edge (11); \path[->] (00) edge (01); \end{tikzpicture} &
\begin{tikzpicture}[>=stealth',semithick,auto,a/.style={circle,fill=white,draw},b/.style={circle,fill=white,draw}]\node[a] (00) at (0,0) {};\node[a] (10) at (0,-1) {};\node [b] (01) at (1,0) {};\node [b] (11) at (1,-1) {}; \path[->] (00) edge (01); \path[->] (10) edge (01); \end{tikzpicture} &
\begin{tikzpicture}[>=stealth',semithick,auto,a/.style={circle,fill=white,draw},b/.style={circle,fill=white,draw}]\node[a] (00) at (0,0) {};\node[a] (10) at (0,-1) {};\node [b] (01) at (1,0) {};\node [b] (11) at (1,-1) {}; \path[->] (00) edge (01); \path[->] (10) edge (11); \end{tikzpicture} &
\begin{tikzpicture}[>=stealth',semithick,auto,a/.style={circle,fill=white,draw},b/.style={circle,fill=white,draw}]\node[a] (00) at (0,0) {};\node[a] (10) at (0,-1) {};\node [b] (01) at (1,0) {};\node [b] (11) at (1,-1) {}; \path[->] (00) edge (11); \path[->] (10) edge (11); \end{tikzpicture}
\\ \hhline{====} & & & \\[0.1em]
\begin{tikzpicture}[>=stealth',semithick,auto,a/.style={circle,fill=white,draw},b/.style={circle,fill=white,draw}]\node[a] (00) at (0,0) {};\node[a] (10) at (0,-1) {};\node [b] (01) at (1,0) {};\node [b] (11) at (1,-1) {}; \path[->] (00) edge (01); \path[->] (10) edge (01); \end{tikzpicture} &
\begin{tikzpicture}[thick,auto,a/.style={circle,fill=white,draw},b/.style={circle,fill=white,draw}]\node[a] (00) at (0,0) {};\node[a] (10) at (0,-1) {};\node [a] (01) at (1,0) {};\node [a] (11) at (1,-1) {};\node[b] (02) at (2,0){};\node[b] (12) at (2,-1){};\draw (00) edge (01);\draw (01) edge (02);\draw (00) edge[bend left] (02);\end{tikzpicture}  &
\begin{tikzpicture}[thick,auto,a/.style={circle,fill=white,draw},b/.style={circle,fill=white,draw}]\node[a] (00) at (0,0) {};\node[a] (10) at (0,-1) {};\node [a] (01) at (1,0) {};\node [a] (11) at (1,-1) {};\node[b] (02) at (2,0){};\node[b] (12) at (2,-1){};\draw (00) edge (01);\draw (01) edge (02);\draw (00) edge[bend left] (02);\end{tikzpicture}  &
\begin{tikzpicture}[thick,auto,a/.style={circle,fill=white,draw},b/.style={circle,fill=white,draw}]\node[a] (00) at (0,0) {};\node[a] (10) at (0,-1) {};\node [a] (01) at (1,0) {};\node [a] (11) at (1,-1) {};\node[b] (02) at (2,0){};\node[b] (12) at (2,-1){};\draw (10) edge (01);\draw (01) edge (12);\draw (10) edge[bend right] (12);\end{tikzpicture}
\\ \hhline{====} & & & \\[0.1em]
\begin{tikzpicture}[>=stealth',semithick,auto,a/.style={circle,fill=white,draw},b/.style={circle,fill=white,draw}]\node[a] (00) at (0,0) {};\node[a] (10) at (0,-1) {};\node [b] (01) at (1,0) {};\node [b] (11) at (1,-1) {}; \path[->] (00) edge (01); \path[->] (10) edge (11); \end{tikzpicture} &
\begin{tikzpicture}[thick,auto,a/.style={circle,fill=white,draw},b/.style={circle,fill=white,draw}]\node[a] (00) at (0,0) {};\node[a] (10) at (0,-1) {};\node [a] (01) at (1,0) {};\node [a] (11) at (1,-1) {};\node[b] (02) at (2,0){};\node[b] (12) at (2,-1){};\draw (00) edge (01);\draw (01) edge (02);\draw (00) edge[bend left] (02);\end{tikzpicture}  &
\begin{tikzpicture}[thick,auto,a/.style={circle,fill=white,draw},b/.style={circle,fill=white,draw}]\node[a] (00) at (0,0) {};\node[a] (10) at (0,-1) {};\node [a] (01) at (1,0) {};\node [a] (11) at (1,-1) {};\node[b] (02) at (2,0){};\node[b] (12) at (2,-1){};\draw (00) edge (01);\draw (01) edge (02);\draw (00) edge[bend left] (02);\draw (10) edge (11);\draw (11) edge (12);\draw (10) edge[bend left] (12);\end{tikzpicture}  &
\begin{tikzpicture}[thick,auto,a/.style={circle,fill=white,draw},b/.style={circle,fill=white,draw}]\node[a] (00) at (0,0) {};\node[a] (10) at (0,-1) {};\node [a] (01) at (1,0) {};\node [a] (11) at (1,-1) {};\node[b] (02) at (2,0){};\node[b] (12) at (2,-1){};\draw (10) edge (11);\draw (11) edge (12);\draw (10) edge[bend left] (12);\end{tikzpicture}
\\ \hhline{====} & & & \\[0.1em]
\begin{tikzpicture}[>=stealth',semithick,auto,a/.style={circle,fill=white,draw},b/.style={circle,fill=white,draw}]\node[a] (00) at (0,0) {};\node[a] (10) at (0,-1) {};\node [b] (01) at (1,0) {};\node [b] (11) at (1,-1) {}; \path[->] (00) edge (11); \path[->] (10) edge (11); \end{tikzpicture} &
\begin{tikzpicture}[thick,auto,a/.style={circle,fill=white,draw},b/.style={circle,fill=white,draw}]\node[a] (00) at (0,0) {};\node[a] (10) at (0,-1) {};\node [a] (01) at (1,0) {};\node [a] (11) at (1,-1) {};\node[b] (02) at (2,0){};\node[b] (12) at (2,-1){};\draw (00) edge (11);\draw (11) edge (02);\draw (00) edge[bend left] (02);\end{tikzpicture}  &
\begin{tikzpicture}[thick,auto,a/.style={circle,fill=white,draw},b/.style={circle,fill=white,draw}]\node[a] (00) at (0,0) {};\node[a] (10) at (0,-1) {};\node [a] (01) at (1,0) {};\node [a] (11) at (1,-1) {};\node[b] (02) at (2,0){};\node[b] (12) at (2,-1){};\draw (10) edge (11);\draw (11) edge (12);\draw (10) edge[bend left] (12);\end{tikzpicture}  &
\begin{tikzpicture}[thick,auto,a/.style={circle,fill=white,draw},b/.style={circle,fill=white,draw}]\node[a] (00) at (0,0) {};\node[a] (10) at (0,-1) {};\node [a] (01) at (1,0) {};\node [a] (11) at (1,-1) {};\node[b] (02) at (2,0){};\node[b] (12) at (2,-1){};\draw (10) edge (11);\draw (11) edge (12);\draw (10) edge[bend left] (12);\end{tikzpicture}
\\ \hhline{----} \\
$A\quad B$
\end{tabular}

\vskip1cm

\begin{tabular}{c||c|c|c}
$C\quad A$
\\ \hhline{----}
\begin{tikzpicture}[>=stealth',semithick,auto,a/.style={circle,fill=white,draw},b/.style={circle,fill=white,draw}]\node[a] (00) at (0,0) {};\node[a] (10) at (0,-1) {};\node [b] (01) at (1,0) {};\node [b] (11) at (1,-1) {}; \path[->] (10) edge (11); \path[->] (00) edge (11); \end{tikzpicture} &
\begin{tikzpicture}[>=stealth',semithick,auto,a/.style={circle,fill=white,draw},b/.style={circle,fill=white,draw}]\node[a] (00) at (0,0) {};\node[a] (10) at (0,-1) {};\node [b] (01) at (1,0) {};\node [b] (11) at (1,-1) {}; \path[->] (00) edge (01); \path[->] (10) edge (01); \end{tikzpicture} &
\begin{tikzpicture}[>=stealth',semithick,auto,a/.style={circle,fill=white,draw},b/.style={circle,fill=white,draw}]\node[a] (00) at (0,0) {};\node[a] (10) at (0,-1) {};\node [b] (01) at (1,0) {};\node [b] (11) at (1,-1) {}; \path[->] (00) edge (01); \path[->] (10) edge (11); \end{tikzpicture} &
\begin{tikzpicture}[>=stealth',semithick,auto,a/.style={circle,fill=white,draw},b/.style={circle,fill=white,draw}]\node[a] (00) at (0,0) {};\node[a] (10) at (0,-1) {};\node [b] (01) at (1,0) {};\node [b] (11) at (1,-1) {}; \path[->] (00) edge (11); \path[->] (10) edge (11); \end{tikzpicture}
\\ \hhline{====} & & & \\[0.1em]
\begin{tikzpicture}[>=stealth',semithick,auto,a/.style={circle,fill=white,draw},b/.style={circle,fill=white,draw}]\node[a] (00) at (0,0) {};\node[a] (10) at (0,-1) {};\node [b] (01) at (1,0) {};\node [b] (11) at (1,-1) {}; \path[->] (00) edge (01); \path[->] (10) edge (01); \end{tikzpicture} &
\begin{tikzpicture}[thick,auto,a/.style={circle,fill=white,draw},b/.style={circle,fill=white,draw}]\node[a] (00) at (0,0) {};\node[a] (10) at (0,-1) {};\node [a] (01) at (1,0) {};\node [a] (11) at (1,-1) {};\node[b] (02) at (2,0){};\node[b] (12) at (2,-1){};\draw (10) edge (01);\draw (01) edge (02);\draw (10) edge (02);\end{tikzpicture}  &
\begin{tikzpicture}[thick,auto,a/.style={circle,fill=white,draw},b/.style={circle,fill=white,draw}]\node[a] (00) at (0,0) {};\node[a] (10) at (0,-1) {};\node [a] (01) at (1,0) {};\node [a] (11) at (1,-1) {};\node[b] (02) at (2,0){};\node[b] (12) at (2,-1){};\draw (10) edge (01);\draw (01) edge (02);\draw (10) edge (02);\end{tikzpicture}  &
\begin{tikzpicture}[thick,auto,a/.style={circle,fill=white,draw},b/.style={circle,fill=white,draw}]\node[a] (00) at (0,0) {};\node[a] (10) at (0,-1) {};\node [a] (01) at (1,0) {};\node [a] (11) at (1,-1) {};\node[b] (02) at (2,0){};\node[b] (12) at (2,-1){};\draw (10) edge (01);\draw (01) edge (12);\draw (10) edge[bend right] (12);\end{tikzpicture}
\\ \hhline{====} & & & \\[0.1em]
\begin{tikzpicture}[>=stealth',semithick,auto,a/.style={circle,fill=white,draw},b/.style={circle,fill=white,draw}]\node[a] (00) at (0,0) {};\node[a] (10) at (0,-1) {};\node [b] (01) at (1,0) {};\node [b] (11) at (1,-1) {}; \path[->] (00) edge (01); \path[->] (10) edge (11); \end{tikzpicture} &
\begin{tikzpicture}[thick,auto,a/.style={circle,fill=white,draw},b/.style={circle,fill=white,draw}]\node[a] (00) at (0,0) {};\node[a] (10) at (0,-1) {};\node [a] (01) at (1,0) {};\node [a] (11) at (1,-1) {};\node[b] (02) at (2,0){};\node[b] (12) at (2,-1){};\draw (10) edge (11);\draw (11) edge (02);\draw (10) edge (02);\end{tikzpicture}  &
\begin{tikzpicture}[thick,auto,a/.style={circle,fill=white,draw},b/.style={circle,fill=white,draw}]\node[a] (00) at (0,0) {};\node[a] (10) at (0,-1) {};\node [a] (01) at (1,0) {};\node [a] (11) at (1,-1) {};\node[b] (02) at (2,0){};\node[b] (12) at (2,-1){};\draw (10) edge (11);\draw (11) edge (12);\draw (10) edge[bend left] (12);\end{tikzpicture}  &
\begin{tikzpicture}[thick,auto,a/.style={circle,fill=white,draw},b/.style={circle,fill=white,draw}]\node[a] (00) at (0,0) {};\node[a] (10) at (0,-1) {};\node [a] (01) at (1,0) {};\node [a] (11) at (1,-1) {};\node[b] (02) at (2,0){};\node[b] (12) at (2,-1){};\draw (10) edge (11);\draw (11) edge (12);\draw (10) edge[bend left] (12);\end{tikzpicture}
\\ \hhline{====} & & & \\[0.1em]
\begin{tikzpicture}[>=stealth',semithick,auto,a/.style={circle,fill=white,draw},b/.style={circle,fill=white,draw}]\node[a] (00) at (0,0) {};\node[a] (10) at (0,-1) {};\node [b] (01) at (1,0) {};\node [b] (11) at (1,-1) {}; \path[->] (00) edge (11); \path[->] (10) edge (11); \end{tikzpicture} &
\begin{tikzpicture}[thick,auto,a/.style={circle,fill=white,draw},b/.style={circle,fill=white,draw}]\node[a] (00) at (0,0) {};\node[a] (10) at (0,-1) {};\node [a] (01) at (1,0) {};\node [a] (11) at (1,-1) {};\node[b] (02) at (2,0){};\node[b] (12) at (2,-1){};\draw (10) edge (11);\draw (11) edge (02);\draw (10) edge (02);\end{tikzpicture}  &
\begin{tikzpicture}[thick,auto,a/.style={circle,fill=white,draw},b/.style={circle,fill=white,draw}]\node[a] (00) at (0,0) {};\node[a] (10) at (0,-1) {};\node [a] (01) at (1,0) {};\node [a] (11) at (1,-1) {};\node[b] (02) at (2,0){};\node[b] (12) at (2,-1){};\draw (10) edge (11);\draw (11) edge (12);\draw (10) edge[bend left] (12);\end{tikzpicture}  &
\begin{tikzpicture}[thick,auto,a/.style={circle,fill=white,draw},b/.style={circle,fill=white,draw}]\node[a] (00) at (0,0) {};\node[a] (10) at (0,-1) {};\node [a] (01) at (1,0) {};\node [a] (11) at (1,-1) {};\node[b] (02) at (2,0){};\node[b] (12) at (2,-1){};\draw (10) edge (11);\draw (11) edge (12);\draw (10) edge[bend left] (12);\end{tikzpicture}
\\ \hhline{----} \\
$A\quad B$
\end{tabular}

\vskip1cm
\[
\begin{matrix}{
\begin{tabular}{|c|c|c|}
\hhline{---} & &\\
$o^B_{111}$ & $o^B_{111}$ & $o^A_{112}$
\\ \hhline{---} & &\\
$o^B_{111}$&$o^B_{111}$&$o^A_{112}$
\\ \hhline{---} & &\\
$o^C_{121}$&$o^B_{122}$&$o^B_{122}$
\\ \hhline{---}
\end{tabular}
}
&{\phantom{abcd}}&
{
\begin{tabular}{|c|c|c|}
\hhline{---} & &\\
$o^A_{111}$ & $o^A_{111}$ & $o^C_{212}$
\\ \hhline{---} & &\\
$o^A_{111}$&$o^A_{111}$&$o^A_{222}$
\\ \hhline{---} & &\\
$o^C_{121}$&$o^C_{222}$&$o^A_{222}$
\\ \hhline{---}
\end{tabular}
}
&{\phantom{abcd}}&
{
\begin{tabular}{|c|c|c|}
\hhline{---} & &\\
$o^B_{211}$ & $o^B_{211}$ & $o^C_{212}$
\\ \hhline{---} & &\\
$o^B_{221}$&$o^C_{222}$&$o^A_{222}$
\\ \hhline{---} & &\\
$o^B_{221}$&$o^C_{222}$&$o^C_{222}$
\\ \hhline{---}
\end{tabular}
}
\end{matrix}
\]

\medskip

Figure 8. A monotone bargaining three person game game correspondence 
which is not tight and not NS, 
 
\bigskip

In this figure 8, an upper index $A, B$, or $C$
over a strategy profile $(x,y,z)$, shows
which player, Alice, Bob, or Claire, could improve $(x,y,z)$
by changing his/her strategy.
Since each profile has an upper index, we conclude  that
the game has no NE
for any realization of the considered partially defined utility function
$u = (u_A, u_B, u_C)$.

\smallskip

Let us note that in the three-person case
tightness and NS are no longer related.
Tightness is neither necessary nor sufficient for NS;
see examples in \cite{Gur75} and \cite{Gur88}.

\begin{remark}
In \cite[Remark 3]{Gur75} it was mistakenly claimed that
tightness is necessary for NS even in the $n$-person case.
This mistake was noticed by Danilov and corrected in \cite{Gur88}.
\end{remark}


{\bf Acknowledgements:} We are thankful to Endre Boros for helpful remarks and suggestions and Denis Mironov for help in drawing Pictures.
{The first author was partially funded by Russian Academic Excellence Project '5-100'.}

\end{document}